\documentclass[a4paper]{article}
\usepackage{latexsym}
\usepackage{amsbsy}
\usepackage{amssymb}
\usepackage{amsmath}
\usepackage{graphicx}
\usepackage{epsfig}
\usepackage[latin1]{inputenc}
\usepackage{geometry}
\usepackage{caption}
\usepackage{epstopdf}
\usepackage{float}
\usepackage{xcolor}

\input{amssym.def}
\input{amssym}

\makeatletter
\@namedef{subjclassname@2020}{\textup{2020} Mathematics Subject Classification}
\makeatother

\newtheorem{thm}{Theorem}[section]
\newtheorem{lemma}[thm]{Lemma}

\newtheorem{definition}[thm]{Definition}

\newtheorem{problem}{Problem}
\newtheorem{remark}[thm]{Remark}
\newcommand{\pr}{\noindent{\bf Proof. }}
\newcommand{\ep}{\nolinebreak{\hspace*{\fill}$\Box$ \vspace*{0.25cm}}\\}

\newcommand{\R}{\mathbb{R}}
\newcommand{\Rn}{\mathbb{R}^n}
\newcommand{\N}{\mathbb{N}}
\newcommand{\C}{\mathbb{C}}
\newcommand{\D}{\mathcal{D}}

\newcommand{\E}{\mathcal{E}}
\newcommand{\eps}{\varepsilon}

\newcommand{\net}[1]{(#1)_{\eps}}

\DeclareMathOperator{\supp}{supp}

\let \Re \relax
\DeclareMathOperator{\Re}{Re}

\usepackage{geometry}

\usepackage{subfigure}

\def\dj{d\kern-0.4em\char"16\kern-0.1em}
\def\Dj{\mbox{\raise0.3ex\hbox{-}\kern-0.4em D}}
\begin{document}

\title{\textbf{The Euler-Bernoulli equation \break with distributional coefficients and forces}}

\author{Robin Blommaert\thanks{
Department of Mathematics: Analysis, Logic and Discrete Mathematics, University of Ghent, Krijgslaan 281 (building S8), 9000 Ghent, Belgium, Electronic mail: Robin.Blommaert@gmail.com},
Sr\dj an Lazendi\'c\thanks{Department of Electronics and Information Systems, Clifford Research Group, University of Ghent, Krijgslaan 281 (building S8), 9000 Ghent, Belgium, Electronic mail: Srdan.Lazendic@UGent.be (Corresponding author)}, and
Ljubica Oparnica\thanks{Department of Mathematics: Analysis, Logic and Discrete Mathematics, University of Ghent, Krijgslaan 281 (building S8), 9000 Ghent, Belgium \& Faculty of Education, University of Novi Sad, Podgori\v cka 4, 25000 Sombor, Serbia, Electronic mail: oparnica.ljubica@ugent.be, ljubica.oparnica@gmail.com} }
\date{}
\maketitle

\begin{abstract}
\noindent In this work we investigate a very weak solution to the initial-boundary value problem of an Euler-Bernoulli beam model. We allow for bending stiffness, axial- and transversal forces as well as for initial conditions to be irregular functions or distributions. We prove the well-posedness of this problem in the very weak sense. More precisely, we define the very weak solution to the problem and show its existence and uniqueness. For regular enough coefficients we show consistency with the weak solution. Numerical analysis shows that the very weak solution coincides with the weak solution, when the latter exists, but also offers more insights into the behaviour of the very weak solution, when the weak solution doesn't exist.
\end{abstract}
	
\tableofcontents
	
\section{Introduction}
In this paper we investigate  the  initial-boundary value problem which arises in the analysis of the Euler-Bernoulli beam equation that includes irregular coefficients and force terms. More precisely, we consider the initial-boundary value problem
\begin{align}\label{P1}
	& \partial_x^2\left(c(x) \partial_x^2 u \right) + b(x,t) \partial_x^2 u +\partial_t^2 u = g,\\
	& u |_{t=0} = f_1, \qquad \partial_t u|_{t=0} = f_2,\label{P1}\\
	& u  |_{x=0}= u |_{x=1} = 0, \qquad  \partial_x u |_{x=0} = \partial_x u |_{x=1} = 0,
\end{align}
	where the coefficients \(c, b\), the force term \(g\) and the initial data \(f_1\) and \(f_2\) are irregular objects, such as Dirac-type distributions. The structure of the considered problem arises from mechanics and it describes the displacement of a beam under axial and transversal forces.
	
	Similar problems, involving distributional and generalized solutions to the Euler-Bernoulli model have already been addressed in literature \cite{BC2007, HO2009, HKO2013, YS2001}. In the Colombeau algebra \cite{GKOS2001} setting, it was shown that if \(c\), \(b\), \( g\)  and \(f_1\) and \(f_2\) are Colombeau generalized functions, the problem (1)-(3) has a unique solution in the Colombeau functions based on families of weak solutions to the corresponding regularized problem \cite{HO2009}. 
	
	Recently introduced concept of the very weak solution \cite{GR15} offers an alternative approach to deal with irregular objects. Since then the very weak solution has been used to successfully treat a wide range of problems \cite{ARST2020c, ARST2021, ARST2021b,Gar20, MRT2019a, MRT2019b, RT17a, RT17b, RT2018, RY2020} with irregular data. The basic idea is to regularize irregular objects  by nets of regular functions with moderate asymptotics, find the weak solution to the net of regularized problems obtaining a net of solutions, which is called a \textit{sequential solution}. When the sequential solution is of moderate growth, it is called a \textit{very weak solution}. 
		
	To the best of our knowledge, the concept of the very weak solution for the Euler-Bernoulli beam initial-boundary value problem with irregular coefficients and forces has not been addressed so far. The results presented in \cite{HO2009} include a weak solution of the problem (1)-(3) for regular data together with the Colombeau generalized solutions. In comparison with the Colombeau algebra generalized solution approach \cite{HO2009}, in the concept of very weak solution, we do not require derivatives to be of moderate growth and the approximating nets are chosen more freely, independently of the fact whether the underlying space is an algebra or not \cite{GR15}.
	
	The outline of the article is as follows. First, we conclude the Introduction section by fixing some notations, recall main notions and definitions important for our work, and give a brief description of the mechanical model for the Euler-Bernoulli beam. In Section 2, we review the existing results on the weak solution for the $L^\infty$-coefficients obtained in \cite{HO2009}, and modify the existing energy estimate making it suitable for the very weak solution approach. In Section 3 we define and show existence of the very weak solution. Section 4 includes the uniqueness results of the very weak solution defined in an appropriate sense and the consistency results of the very weak solution with the weak solution, when the coefficients are regular enough. In Section 5 we conduct numerical analysis of the obtained theoretical results. Finally, we give a short conclusion and discuss possibilities for future research.

	\subsection{Notations and notions} 	Let \(L^\infty(\R^d)\) be the space of essentially bounded measurable functions, \(C^k(\R^d)\) is the space of $k$-times continuously differentiable functions on $\mathbb{R}^d$, \(\D(\R^d) := C^\infty_0(\R^d)\) the space of compactly supported smooth functions, \(\D'(\R^d)\) the space of distributions, and \(\E'(\R^d)\) the space of compactly supported distributions. Further let \(L^2(\R^d)\) be the space of square integrable functions, \(H^2(\R^d)\) be the Sobolev space that consists of the functions \(v \in L^2(\R^d)\) such that its weak derivatives up to order 2 exist and are square integrable functions. \(H^2_0(\R^d)\) denotes the closure of \(\D(\R^d)\) in \(H^2(\R^d)\). The space \(H^{-2}(\R^d)\) is the space of continuous linear functionals \(f: H^2(\R^d)\to \C\) with the norm
	\begin{align*}
	\|f\|_{H^{-2}(\R^d)}:= \sup\limits_{v \in H^2(\R^d)} \dfrac{|f(v)|} {\|v\|_{H^2(\R^d)}}.
	\end{align*}
	There are vector valued Sobolev spaces that are essential in the construction of weak and very weak solutions to PDEs, when working with time dependent problems. Let \((X, \|\cdot\|_X)\) be a Banach space, then \(L^2([0,T], X)\) is the space of functions \(f\) over \([0,T]\) with values in \(X\) with the norm 
	\begin{align*}
	\|f\|_{L^2([0,T], X)} := \left( \int_0^T \|f(t)\|_X dt \right)^{1/2}<\infty.
	\end{align*}
	The space \(H^2([0,T], X)\) is the subspace consisting of all functions \(f\) in \(L^2([0,T],X)\) for which the norm
	\begin{align*}
	\|f\|_{H^2([0,T], X)}:= \left( \int_0^T \|f(t)\|_X dt + \int_0^T \|\partial_t f(t)\|_X dt +  \int_0^T \|\partial_t^2 f(t)\|_X dt \right)^{1/2}
	\end{align*}
	is finite, where \(\partial_t\) is interpreted as the weak derivative.
	Similarly \(C([0,T], X)\) is the space of continuous functions over \([0,T]\) with values in \(X\) endowed with the norm
	\begin{align*}
	\|f\|_{C([0,T], X)} := \sup_{t \in [0,T]} \|f(t)\|_X<\infty.
	\end{align*}
 	For our purposes, the Hilbert space  \(L^2(0,1)\) with the inner product given by $$\langle u, v\rangle := \int_0^1 u(x)\overline{v(x)} \,dx$$ will be of special interest. Together with this inner product we consider the anti-dual action of \(f \in H^{-2}(0,1)\) by 	$\langle f, v\rangle := f(\overline{v}),$ for all $v \in H^2(0,1)$.
 
For the sake of simplifying certain notations, we will also use the Landau asymptotic symbols. More precisely, we say that function $f$ is ''big-O'' of function $g$,  denoted by $f(x)=O(g(x))$, when $x\to 0$, if there exist positive numbers $\delta$ and $M$ such that $$\left\vert f(x)\right\vert\leq M\left\vert g(x)\right\vert,~~~x\in(-\delta,\delta)\backslash\{0\}.$$ Similarly, if $g(x)\not=0$, for $x\not=0$, we say that function $f$ is ''small-o'' of function $g$, as $x\to 0$, if $$\lim_{x\to 0} \frac{\left\vert f(x)\right\vert}{\left\vert g(x)\right\vert}=0.$$

		\subsection{The Euler-Bernoulli beam model}\label{secModel}
	In this subsection we will explain in more detail the model of an elastic rod behind the Euler-Bernoulli beam equation. We base ourselves on a general theory of elasticity, which can be found in \cite{a-ster,AG2000}.
	
	Let us consider a cylindrical beam of length \(L\). A beam is often modeled by considering deviations from the beam axis. Let us assume that the beam is under influence of a vertical force \(g_1\) and an axial force \(P\), but there are no lateral forces. The vertical displacement of the beam axis in this situation, depends on time \(t\) and position \(x\) and will be modeled by the function \(u(x,t)\).
	
	Material particles also experience displacements in the axial direction \(w(x,z,t)\). Let \(p\) be a particle that starts at the position \((x,z)\) at time \(t = 0\). Then \(w(x,z,t)\) is the horizontal displacement of \(p\) at time \(t\). For small horizontal deformations we can approximate the horizontal deformation by \(w(x,z,t) = z \partial_x u(x,t)\). When the elasticity of the beam is modeled by Hooke's law, the stress \(\sigma(x, t)\) is then proportional to the relative extension of the material. Locally, the relative extension is given by the derivative \(\partial_x w(x,z,t)\). This together implies that $$\sigma(x,t) = E(x) \partial_x w(x,z,t),$$ where the proportionality constant \(E(x)\) is called the modulus of elasticity. The bending moment \(M\) is a stress resultant given by 
	\(M(x,t) = \int_{C(x)}\sigma(x,t) z dA\), with the integral over the cross-section of the beam \(C(x)\) at position \(x\). In our case, we obtain
	\begin{align*}
	M(x,t) &= \int_{C(x)} E(x) \partial_x w(x,z,t) z dA = \int_{C(x)} E(x) \partial_x^2 u(x,t) z^2 dA \\
	&= E(x) \partial_x^2 u(x,t) \int_{C(x)} z^2 dy dz = E(x) I(x) \partial_x^2 u(x,t).
	\end{align*}
	Note that \(I(x):= \int_{C(x)} z^2 dy dz\) is the second moment of area of the beam's cross-section. A shearing force is defined as an excess horizontal force on the top part of the beam over the bottom part of the beam. The bending moment defines the shearing force as
	\[
	Q_1(x,t) = -\partial^2_x M(x,t) = -\partial^2_x\left(E(x) I(x) \partial_x^2 u(x,t)\right).
	\]
	The minus sign appears as an excess force in a position \(x\) with positive curvature results in a downwards force.
	Next we have the shearing force due to the axial force, that is
	\[
	Q_2(x,t) = - P(t)\partial^2_x u(x,t).
	\]
	Using Newton's second law on the vertical displacement \(u(x,t)\) and the vertical forces \(g_1(x,t)\), \(Q_1(x,t)\) and \(Q_2(x,t)\), results in the dynamic Euler-Bernoulli equation
	\begin{equation} \label{EBNewton}
	R(x)\partial_t^2 u = g_1(x,t) + Q_2(x,t) + Q_2(x,t).
	\end{equation}
	The line density \(R(x)\) is a replacement for the proper mass, as we are working locally on a cross-section of infinitesimal width. It is defined as \(\partial_x m(x)\) with \(m(x)\) the mass of the beam contained in the interval \([0,x]\). By putting together the previous results, we obtain the dynamic Euler-Bernoulli equation of the form 
\begin{equation*}
	\partial_x^2\left(A(x) \partial_x^2 u\right) + P(t) \partial_x^2 u + R(x)\partial_t^2 u = g_1(x,t),\quad x\in (0,1),\quad t>0, 
\end{equation*} where $A(x):=E(x)I(x)$ is the bending stiffness, $P(t)$ is the axial force, $R(x)$ is the line density and $g_1(x,t)$ denotes the force term.

	In this article we will consider the Euler-Bernoulli beam model where now the rod consists of two parts with uniform cross-section. Moreover, we consider \textit{distributional} trans\-versal force \(g_1\) and axial force \(P\). Denote the vertical displacement of the beam axis at \((x,t)\) by \(u(x,t)\). The equation of motion then reads as follows
	\begin{equation}
	\partial_x^2\left(A(x) \partial_x^2 u\right) + P(t) \partial_x^2 u + R(x)\partial_t^2 u = g_1(x,t),\quad x\in (0,1),\quad t>0, \label{eqEBModel}
	\end{equation}
	where by taking into account that the beam consists of two parts, the following holds:
	\begin{itemize}
	\item \(A\) is the bending stiffness, given by $A(x)=EI_1+H(x-x_0)(EI_2-EI_1)$ where $H$ is the Heaviside jump function and $I_1\not= I_2$ are are the second moments of area of the beam's cross-section that correspond to the two parts of the beam,
	\item \(R\) is the line density of the beam of the form $R(x)=R_0+H(x-x_0)(R_2-R_1)$,
	\item \(P\) is an axial force of the form $P(t)=P_0+P_1\delta(t-t_1)$, $P_0,P_1>0$, 
	\item \(g_1\) is the vertical force term.
	\end{itemize}
	Further, we consider the \textit{clamped Euler-Bernoulli beam}, i.e., the case when the equation \eqref{eqEBModel} is accompanied by the homogeneous boundary conditions
	\begin{equation*}
	u(0,t) = u(1,t) = \partial_x u(0,t) = \partial_x u(1,t) = 0. \label{eqbc}
	\end{equation*}
	Using a formal change of variables \(t\mapsto \sqrt{R(x)}t\), we reduce the model \eqref{eqEBModel} to the following initial-boundary value problem on the domain \(X_T :=(0,T) \times (0,1)\).
    \begin{problem}\label{problemEB}
     Find a function \(u\) such that
    \begin{equation}
	    \partial_x^2\left(c(x) \partial_x^2 u(x,t) \right) + b(x,t) \partial_x^2 u(x,t) +\partial_t^2 u(x,t) = g(x,t), \quad (x,t) \in X_T,\label{eqEB}\tag{EB}
	\end{equation}
	satisfying the initial conditions
	\begin{equation}
	u(x,0) = f_1(x), \quad \partial_t u(x,0) = f_2(x), \quad x \in [0,1], 
    \label{eqEBic}\tag{IC}
	\end{equation}
	and the boundary conditions
	\begin{equation}
	    u(0,t) = \partial_x u(0,t) = 0,\quad u(1,t) = \partial_x u(1,t) = 0. \label{eqBC} \tag{BC}
	\end{equation}
    \end{problem}
	Note that the function $c$ in \eqref{eqEB} equals $A$, the function $b$ is given by $b(x,t)=P(\sqrt{R(x)}t)$ and \(g(x,t) = g_1(x,\sqrt{R(x)}t)\).
   
In Section \ref{secVWS} we will consider Problem \ref{problemEB} with \(c\in \D'(0,1)\),  \(b\in \D'(X_T)\), \(g\in \D'(X_T)\), $f_1,f_2\in \D'(0,1)$. First, we will say a bit more about about regularization techniques, which are a common approach when dealing with irregular coefficients. Moreover, we will introduce the concept of the very weak solution.
	
	\subsection{A note on regularizations}
	\label{sec:regular}
	A common approach in the functional analysis to deal with irregular objects is via regularization, which will provide us with a regularizing net of smooth functions. We will work with nets of functions \( (u_\eps)_{\eps}\) in a function space \( X \), which is any map \( (0,\eps_0] \to X: \eps \mapsto u_\eps\) for some \(\eps_0>0\). In this article we are concerned by nets of smooth functions indexed by a parameter $\varepsilon$ and obeying certain asymptotic growth estimates as $\varepsilon$ tends to zero.
	
	One can regularize a distribution \(u \in \D'(\R^d)\) by means of a mollifier \(\varphi\). A \textit{Friedrichs mollifier} is a smooth function \(\varphi \in C^{\infty}_0(\R^d)\), \(\varphi \geq 0\) and \(\int \varphi = 1\). The corresponding mollifying net \(\net{\varphi_\eps}\) is 
	\begin{align*}
	\varphi_\eps(x) = \eps^{-d}\varphi(\eps^{-1}x), \qquad \eps>0.
	\end{align*}
	Convolution with the mollifying net produces a regularizing net \(\net{u_\eps}\) of smooth functions by
	\begin{align}
	\label{Mol}
	u_\eps = u*\varphi_\eps, \quad \eps > 0,
	\end{align}
	such that \(u_\eps\to u\) in \(\D'(\R^d)\) as \(\eps\to 0\).
	Notably, the regularizing nets admit asymptotic estimates in terms of the regularization parameter. 
	\begin{definition}[Moderateness]\label{Mod}
		Let \((X, \|\cdot\|_X)\) be a Banach space. A net of elements \(\net{v_\eps}\) in \(X\) is \(X\)-moderate if there exist \(C> 0\) and $N\in\mathbb{N}_{0}$ such that
		\begin{align*}
		\|v_\eps\|_X  \leq C \eps^{-N}, \quad \text{as } \eps \to 0.
		\end{align*}
	\end{definition}
	
\begin{remark} \label{reModerateMollification}
Specifically, the net \(\net{u_\eps}\) defined by \eqref{Mol} and all nets of partial derivatives \(\net{\partial^\alpha_x u_\eps}\) are \(L^p\)-moderate on compact sets, for \(1\leq p\leq \infty\).
Indeed, let \(v \in \E'(\R^d)\). By the structure theorem for compactly supported distributions of finite order \cite{GKOS2001,V1979}, we can write \(v\) as a finite sum of distributional derivatives
		\[
		v= \sum_{\alpha \leq m}\partial^\alpha f,
		\]
		on \(\supp v\), for some \(m \in \N\) and a function \(f \in C(\Rn)\) with compact support in an arbitrarily small neighborhood of \(\supp v\).
		Convolution with a Friedrichs mollifier then gives
		\begin{align*}
		v_\eps = v*\varphi_\eps &= \sum_{|\alpha| \leq m}\partial^\alpha f*\varphi_\eps\\
		&= \sum_{|\alpha| \leq m} f*\partial^\alpha\varphi_\eps \\
		&= \sum_{|\alpha| \leq m}\eps^{-|\alpha|} f*(\eps^{-d}(\partial^\alpha\varphi)(x/\eps)),
		\end{align*}
		such that by Young's inequality
		\begin{align*}
		    \|v_\eps\|_{L^p(\R^d)}&\leq \sum_{|\alpha|\leq m} \eps^{-|\alpha|}\|f\|_{L^1(\R^d)}\|\eps^{-d}(\partial^\alpha\varphi)(x/\eps)\|_{L^p(\R^d)}\\
		    &\leq \eps^{-m}\|f\|_{L^1(\R^d)} \sum_{|\alpha| \leq m}\|\partial^\alpha\varphi(x)\|_{L^p(\R^d)}.
		\end{align*}
    Similarly, nets of any order of derivatives \(\net{\partial^\beta v_\eps}\) are \(L^p(\R^d)\)-moderate. Now let us consider any compact \(K\subset \R^d\) and a distribution \(u \in \D'(\R^d)\). By writing \(u\) as a locally finite sum of compactly supported distributions one finds that the regularizing net \(\net{u_\eps}\) is \(L^p(K)\)-moderate. We can conclude that nets with moderate asymptotics are well suited for representing distributional objects.
	\end{remark}

	 In our investigation of the very weak solution we will encounter also log-type moderate nets. Due to its importance in the sequel, we will also show that every distribution of finite order can be regularized such that
the obtained regularization net is of log-type.
	 
	\begin{definition}[Moderateness of log-type]
		Let \((X, \|\cdot \|_X)\) be a Banach space. A net \(\net{v_\eps}\) is \(X\)-moderate of log-type if there is \(C > 0\) such that
		\[
		\|v_\eps\|_X \leq C \log\dfrac{1}{\eps}, \quad \text{as }\eps \to 0.
		\]
	\end{definition}
	
\begin{remark}  It has been already shown \cite{MO1989} that log-type regularizations of distributions can be easily obtained by convolution with logarithmically scaled Friedrichs mollifiers. For the convenience of the reader we will include an explicit construction starting from a moderate net.

Let \((X, \|\cdot \|_X)\) be a Banach space. Suppose that \(\net{a_\varepsilon}\) is an \(X\)-moderate net. By moderateness of \(\net{a_\varepsilon}\) we have the estimate
	\begin{equation}
	\|a_\varepsilon\|_X \leq C_1 \varepsilon^{-N}, \quad \text{as }\eps \to 0,\label{eq: logmoderate1}
	\end{equation}
	for $C_1>0$ and $N\in\mathbb{N}_0$. We will now change the parametrisation of \(\varepsilon\) to get a slower rate of growth. Let us consider the net given by
	\begin{equation}
	\net{b_\varepsilon} = \net{a_{\lambda_\varepsilon}}, \quad \text{where}\quad  \lambda_\varepsilon = \left(\log \dfrac{1}{\varepsilon}	\right)^{-\dfrac{1}{N}},\quad \eps >0.\label{eq: logmoderate8}
	\end{equation}
	It is then clear that \(\net{b_\varepsilon}\) is moderate of log-type. Indeed, by \eqref{eq: logmoderate1} and \eqref{eq: logmoderate8} we obtain that when $\eps\to 0$, the following estimate holds:
	\begin{align*}
	\|b_\varepsilon\|_X = \|a_{\lambda_\varepsilon}\|_X \leq C_1 \lambda_\varepsilon^{-N} = C_1 \left(\left(\log\dfrac{1}{\varepsilon}\right)^{-\dfrac{1}{N}}\right)^{-N} = C_1 \log \dfrac{1}{\varepsilon}.
	\end{align*}
	If the constant \(N\) is unknown, then we can use the following reparametrisation instead. Let 
	\[
	\net{c_\varepsilon} = \net{a_{\mu_\varepsilon}}, \quad \text{with} \quad 	\mu_\varepsilon = \dfrac{1}{\log \log \dfrac{1}{\varepsilon}},\quad \eps >0.
	\]
By using the fact that
	\[
	\left(\log x\right)^N \leq x,
	\]
	for \(x>0\) sufficiently large, we obtain that \(\net{c_\varepsilon}\) is moderate of log-type, since for any \(N\in\mathbb{N}_0\) there holds
	\begin{align*}
	\|c_\varepsilon \|_X = \|a_{\mu_\varepsilon}\|_X \leq C_1\left(\dfrac{1}{\log \log \dfrac{1}{\varepsilon}}\right)^{-N} \leq C_1\left(\log \log \dfrac{1}{\varepsilon}\right)^{N} \leq C_1 \log \dfrac{1}{\varepsilon},
	\end{align*}
	for \(\varepsilon\) sufficiently small. This concludes the proof and the remark.
\end{remark}

\begin{definition}[Negligibility]
Let \((X, \|\cdot\|_X)\) be a Banach space. A net \(\net{v_\eps}\) of elements  in \(X\) is \(X\)-negligible if for all \(q>0\) there exist \(C_q>0\) such that 
\begin{align*}
\|v_\eps\|_X  \leq C_q \eps^{q}, \quad \text{as } \eps \to 0.
\end{align*}
\end{definition}

After introducing the notion of moderateness and neglibility, we have all the ingredients needed to formulate the very weak solution concept.

\subsection{The very weak solution concept}

Now we will explain the general idea behind the very weak solution concept. Although the introduced concepts of moderate and negligible nets are related to those in the Colombeau setting, it is worth noting that the very weak solution concept offers an alternative approach to rigorously analyze the sequential solutions to PDEs. It is also important to note that in the concept of very weak solution, there is no requirement for derivatives to be of moderate growth and the approximating nets can be chosen more freely. 
\begin{itemize}
\item The existence theorem generally has the following structure. Assume that we consider distributional data in the original initial-boundary value problem. Then, we prove that the sequential solution that corresponds to any moderate regularization of the coefficients and data that satisfy some minimal additional conditions, in our case uniform positivity of \(\net{c_\eps}\) and log-type moderateness of \(\net{b_\eps}\), is necessarily moderate. One has to investigate for which classes of irregular functions such regularizations with additional conditions exist and how to construct them. 
\item The uniqueness theorem reads as follows. Let there be two very weak solutions for which the regularizing nets of coefficients and data satisfy the conditions from the existence theorem. The goal is go give sufficient conditions on when the two very weak solutions can be considered the same.
Usually, one shows that when corresponding regularizing nets of the coefficients and data have negligible differences, then the very weak solutions must have negligible differences as well, i.e., the solutions are unique in the negligibility sense. 
\item The main idea of the consistency theorem can be summarized as follows.
The question is whether the very weak solution obtained in the existence theorem, for the case when the data are nice enough functions, coincides with the solution guaranteed by the classical or weak solution theorems. One seeks regularity conditions for the functions or for their regularizing nets such that the very weak solution converges to the known solution.
\end{itemize}

In the following section, we will first discuss the existence and uniqueness of an abstract variational problem and derive energy estimates that will play a fundamental role in the formal introduction of the concept of the very weak solution. After defining the very weak solution to a problem in question, we prove the existence, uniqueness and consistency theorems.

	\section{Abstract variational and Weak solution existence \\ theorems}
	In this section we consider the problem \eqref{eqEB}-\eqref{eqEBic}-\eqref{eqBC} with regular coefficients on an abstract level. More precisely, we give the conditions under which the abstract variational problem as well as the problem \eqref{eqEB}-\eqref{eqEBic}-\eqref{eqBC} has a unique weak solution. Note that Theorem \ref{thmAbstractWeakSolution} is an extension of the abstract variational result from \cite[Thm 1.2, Prop 1.3]{HO2009}, where we consider a more general force term.
	
	\subsection{Abstract variational theorem}
	The abstract variational results presented in \cite{HO2009} are a special case of the more general results from \cite[ch. XVIII,\textsection 5]{DL2000} adapted to the Euler-Bernoulli equation problem. We present a slight adaptation of the result from \cite{HO2009}, which is more suited for our analysis.
	
	Let $V$ and $H$ be two complex, separable Hilbert spaces, where \(V\) is densely embedded into \(H\). 
	If \(V'\) is the anti-dual of \(V\), then \(V\hookrightarrow  H \hookrightarrow  V'\) is a Gelfand triple. 
	Let \(a(t,.,.),\, a_0(t,.,.),\) and \(a_1(t,.,.),\, t\in [0,T],\) be families of continuous sesquilinear 
	forms on \(V\) with
	\begin{equation*}
	a(t,u,v) = a_0(t,u,v) + a_1(t,u,v), \quad \forall u,v \in V,
	\end{equation*}
	such that \(a_0\) and \(a_1\) satisfy the following conditions:
	\begin{enumerate}
		\item[(i)] for all \(u,v\in V: t\mapsto a_0(t,u,v)\) is continuously differentiable \([0,T] \to \C\), i.e., there exist nonnegative constants \(C, C_0\) such that
		\[
		|a_0(t,u,v)|\leq C\|u\|_V \|v\|_V, \quad \text{and}\quad \left|a_0'(t,u,v)\right|\leq  C_0\|u\|_V \|v\|_V,
		\]
		\item[(ii)] \(a_0\) is Hermitian, i.e., \(a_0(t,u,v) = \overline{a_0(t,v,u)}\) for all \(u,v\in V\),
		\item[(iii)] there exist real constants \(\lambda\geq 0\) and \(\alpha>0\) such that 
		\begin{equation*}
		a_0(t,u,u) \geq \alpha \|u\|_V^2 - \lambda \|u\|_H^2, \quad \forall u \in V,\quad  \forall t \in [0,T],
		\end{equation*}
		\item[(iv)] for all \(u, v \in V: t\mapsto a_1(t,u,v)\) is continuous in \([0,T]\to \C\),
		\item[(v)] there exists \(C_1 \geq 0\) such that $$|a_1(t,u, v)| \leq C_1 \|u\|_V \|v\|_H,\quad \forall u,v\in V, \quad \forall t\in[0,T].$$
	\end{enumerate}
	
	\begin{thm}\label{thmAbstractWeakSolution}
		Let \(a(t,.,.),~t\in[0,T]\) satisfy conditions \((i)-(v)\). Let \(u_0\in V,\, u_1 \in H\), and \(f \in L^2([0,T],V')\). Then there exists a unique solution \(u \in L^2([0,T], V)\) satisfying the regularity conditions 
		\begin{equation*}
		\partial_t u \in L^2([0,T], V), \quad \text{and}\quad \partial_t^2 u \in L^2([0,T], V'),
		\end{equation*}
		and solving the abstract initial problem 
		\begin{align}
		&\langle u''(t),v\rangle  + a(t,u(t), v) = \langle f(t), v \rangle, \quad \forall v \in V,\quad \forall t \in (0,T),\\\label{eqAbstractWeakFormulation}
		& u(0)=u_0, \quad u'(0) = u_1.
		\end{align}
		Additionally, the following energy estimate holds
		\begin{equation}
		\|u(t)\|_V^2 + \|u'(t)\|_H^2 \leq \left(D_T \|u_0\|_V^2 + \|u_1\|^2_H  + \int_0^t\|f(\tau)\|^2_{V'} d\tau \right) \cdot \exp(tF_T), \quad \forall t \in[0,T], \label{eqAbstractEstimate}
		\end{equation}
		where the constants \(D_T\) and \(F_T\) are given by
		\begin{gather}
		D_T = (C + \lambda(1+T)) /\min(1,\alpha), \quad \text{and}\label{eqAbstractDT}\\
		F_T = \max(C_0 + C_1, C_1 + 1 + \lambda(1+T))/\min(1,\alpha).\label{eqAbstractFT}
		\end{gather}
	\end{thm}
	
	\begin{remark}
	The differences between Theorem \ref{thmAbstractWeakSolution}
	and the results from \cite{HO2009} are the following.
    \begin{itemize}
	\item[a)] Theorem \ref{thmAbstractWeakSolution} considers the righthandside \(f\) in the space \(L^2([0,T], V')\) whereas \cite{HO2009} considers \(f \in L^2([0,T], H)\). This is justified by \cite[ch. XVIII, \textsection 5, Remark 4.2]{DL2000}. In our case we can use the following estimate
	\[
	2\Re\langle f(\tau) ,u'(\tau)\rangle \leq 2|\langle f(\tau), u'(\tau)\rangle | \leq 2\|f(\tau)\|_{V'}\|u'(\tau)\|_V \leq \|f(\tau)\|^2_{V'} + \|u'(\tau)\|^2_V.
	\]
	to adapt the proof of \cite[Prop. 1.3]{HO2009} and obtain \eqref{eqAbstractEstimate} as an priori estimate. Next, the method of Galerkin approximation given in \cite[ch. XVIII, \textsection 5, Section 2]{DL2000} is applicable and produces a solution with the same properties.
    \item[b)] Compared to \cite{HO2009}, where the constant \(F_T\) is given by \[
    F_T = \max(C_0 + C_1, C_1 + T + 2)/\min(1,\alpha),
    \]
    in our case, the constant $F_T$ in \eqref{eqAbstractFT}, has an additional factor \(\lambda\).
	\end{itemize}
	\end{remark}

	\subsection{The weak solution}\label{secWeakSolution}
In the sequel we will apply the above abstract result to Problem \ref{problemEB} in the case when the coefficients appearing in the equation are regular enough. Thus, let us consider \(H:= L^2(0,1)\) with the standard scalar product and a Sobolev space \(V:= H^2_0(0,1)\). Then, \(V'= H^{-2}(0,1)\) and \(V\hookrightarrow  H \hookrightarrow  V'\) forms a Gelfand triple.	
	
	Formally, taking inner product in \eqref{eqEB} with $v\in H^2_0(0,1)$ and using \eqref{eqBC} we obtain the following weak formulation of the original problem. 
\begin{problem}\label{problemEBWeak}
Find a function $u$ such that 
\begin{align}\label{eqEBWeakFormulation}
	\langle \partial_t^2 u(t),v\rangle + \langle c\partial_x^2 u(t), \partial_x^2 v\rangle + \langle b(t)\partial_x^2 u(t),v \rangle = \langle g(t),v \rangle, \quad t \in (0,T),
\end{align}
satisfying for \(f_1 \in H^2_0(0,1)\) and \(f_2\in L^2(0,1)\)
\begin{align}
u(x,0) = f_1(x), \quad \partial_t u(x,0) = f_2(x),\quad x \in [0,1].\label{eqEBWeakIc}
\end{align}
\end{problem}

\begin{remark} Note that if \(u(t) \in H^2_0(0,1)\), for all $t\in[0,T]$, then the boundary conditions \eqref{eqBC} are automatically satisfied, as we know that $H^2_0(0,1)$ is continuously embedded in the space \\ $\{v\in C^1[0,1]\mid v(0,t)=v(1,t)=\partial_x v(0,t)=\partial_x v(1,t)=0\}$.
\end{remark}

We have now obtained all the ingredients needed to prove the existence and the uniqueness of the weak solution of Problem \ref{problemEB}. In particular, we are able to apply the abstract result from Theorem \ref{thmAbstractWeakSolution} to the initial-value problem specified in Problem \ref{problemEBWeak}.

Thus, for the coefficients \(c \in L^\infty(0,1)\) and \(b \in C([0,T], L^\infty(0,1))\), let us define the sesquilinear forms \(a_0(t,\cdot, \cdot)\) and \(a_1(t,\cdot, \cdot)\) on \(V \times V\) by
	\begin{equation}\label{forms}
	a_0(t,u,v):= \langle c(x)\partial_x^2 u, \partial_x^2 v\rangle, \quad a_1(t,u,v):= \langle b(x,t)\partial_x^2 u, v \rangle, \quad u,v \in H^2_0(0,1),
	\end{equation}
and  \(a(t,\cdot, \cdot) \) with \(a(t,u,v) := a_0(t,u,v) + a_1(t,u,v)\).

\begin{thm}[Weak solution]\label{WeakSolution} 
Let \(c\in L^\infty(0,1)\), \(b\in C([0,T], L^\infty(0,1))\) and \(g\in L^2([0,T], H^{-2}(0,1))\) and suppose there exists a positive constant \(c_0 \in \R\) such that  \[ 0<c_0\leq c(x), \quad \text{for } x\in (0,1).\]		
Let \(f_1\in H^2_0(0,1), f_2 \in L^2(0,1)\). 
Then there exists a unique \(u \in L^2([0,T], H^2_0(0,1))\) with
\begin{equation}
	\partial_t u\in L^2([0,T], H^2_0(0,1)),\quad \partial_t^2 u \in L^2([0,T], H^{-2}(0,1)),\label{eqWeakThmSolSpace}
\end{equation}
	that satisfies \eqref{eqEBWeakFormulation} and initial conditions \(u(0) = f_1\) and $\partial_t u(0) = f_2$. 
	Additionally, for any \(t\in [0,T]\) the energy estimate
\begin{equation}
		\|u(t)\|^2_{H^2(0,1)} + \|\partial_t u(t)\|^2_{L^2(0,1)} \leq \left( D_T \|f_1\|^2_{H^2(0,1)} + \|f_2\|^2_{L^2(0,1)} + \int_{0}^{t} \|g(\tau)\|^2_{H^{-2}(0,1)}d\tau\right) 	\exp(t F_T),\label{EnergyEstimate}
\end{equation}
		holds with the constants
		\begin{align*}
		D_T  = \left(\|c\|_{L^\infty(0,1)} + c_0 C_{1/2} (1+T)\right)/\min\left(\dfrac{c_0}{2},1\right),\\
		F_T = \left(\|b \|_{L^\infty(X_T)} + 1 + c_0 C_{1/2}(1+T)\right)/\min\left(\dfrac{c_0}{2},1\right).
		\end{align*}
	\end{thm}
The proof is a direct consequence of the Theorem \ref{thmAbstractWeakSolution} provided that sesquilinear  forms \eqref{forms} satisfy  properties (i)-(v). We refer to \cite[Theorem 2.2]{HO2009} for the explicit verification of those. The precise abstract constants in the energy estimates read as	
\[
		C  = \|c\|_{L^\infty(0,1)}, \quad C_0 = 0, \quad C_1  =\|b\|_{L^\infty(X_T)},\quad
		\alpha = \dfrac{c_0}{2}, \quad \lambda = c_0C_{1/2},
		\]
where \(C_{1/2}>0\) is an Ehrling's constant, which is defined as the smallest number satisfying
	\[
	\|v\|^2_{H^1(0,1)} \leq \frac{1}{2} \|v\|^2_{H^2(0,1)} + C_{1/2}\|v\|^2_{L^2(0,1)}, \quad v\in H^2_0(0,1).
	\]

From \(u \in L^2([0,T],H^2_0(0,1))\) and \(\partial_t u\in L^2([0,T], H^2_0(0,1))\) and Sobolev embedding theorems we have that  \(u \in C([0,T], H^2_0(0,1))\), while \eqref{eqWeakThmSolSpace} together with \cite[ch. XVIII, \textsection 1, Theorem 1,  pg. 473]{DL2000} implies that \(\partial_t u\in C([0,T],L^2(0,1))\). 
Thus, we have that the solution $u$ belongs to the Banach space denoted by \[W:=C([0,T], H^2_0(0,1))\cap C^1([0,T],L^2(0,1)),\]
endowed with the norm
	\[
	\| u \|_W = \left(\| u \|^2_{L^\infty([0,T],H^2(0,1))} + \|\partial_t u\|^2_{L^\infty([0,T],L^2(0,1))} \right)^{1/2}.
	\]
Finally, note that by taking the supremum over \(t \in [0,T]\) in the energy estimate inequality \eqref{EnergyEstimate}, we easily obtain that
\begin{equation} \label{EnergyEstimate2} 
\|u\|^2_W \leq \left( D_T \|f_1\|^2_{H^2(0,1)} + \|f_2\|^2_{L^2(0,1)} + \|g\|^2_{L^2([0,T], H^{-2}(0,1))} \right) \exp\left(T F_T\right).
\end{equation}

\section{The very weak solution}
	\label{secVWS}
	Now we turn our attention to the case when we have irregular data. 
	Let the coefficients and right hand side in \eqref{eqEB} be distributional, i.e., 
	\begin{equation}\label{eqDistData} b,g \in \D'(X_T),\quad c \in \D'(0,1),\end{equation} 
	and initial data 
	\begin{equation*}
	    f_1 \in \D'(0,1), f_2 \in \D'(0,1).
	\end{equation*}
We introduce the definition of the very weak solution.
		\begin{definition}[Very weak solution]\label{VeryWeakSolutionDef}
		A very weak solution to Problem \ref{problemEB} is a \(W\)-moderate net \(\net{u_\eps}\) if there exist regularizing nets \(\net{f_{1,\eps}}, \net{f_{2,\eps}}\) and \(\net{b_\eps}, \net{g_\eps}\) and \(\net{c_\eps}\) of \(f_1, f_2\), \(b,g\) and \(c\) resp\-ectively such that for every sufficiently small \(\eps>0\),
		the function \(u_\eps\)
		satisfies the regularized gove\-rning equation \begin{equation}\label{eqEBRegWeak} \langle c_\eps(x) \partial_x^2 u_\eps(x,t),\partial_x^2 v(x)\rangle  + \langle b_\eps(x,t) \partial_x^2 u_\eps(x,t), v(x)\rangle + \langle \partial_t^2 u_\eps(x,t), v(x)\rangle  = \langle g_\eps(x,t),v(x)\rangle, \end{equation} for all \(t \in [0,T]\), for all \(v \in H^2_0(0,1)\), and the initial conditions \begin{align}\label{eqEBRegIc} u_\eps(x,0) = f_{1,\eps}(x), \quad \partial_t u_\eps(x,0) = f_{2,\eps}(x),\quad \forall x \in (0,1). \end{align}
	\end{definition} 
	Note that the choice of regularizations of the coefficients \(b\) and \(c\) and right-hand side \(g\) are still unspecified.
    Next, we show that the appropriate moderateness of the coefficients nets yields the existence of a \(W\)-moderate solution net. It then remains to construct such nets as regularization of the distributional data. 
	\begin{thm}[Existence]\label{veryWeakSolutionExists}
		Let \(f_1\), \(f_2\) in \(\D'(0,1)\) and let \(b\), \(c\) and \(g\) satisfy \eqref{eqDistData}.
		\sloppy Let \(\net{f_{1,\eps}}\) be \(H^2_0(0,1)\)-moderate and \(\net{f_{2,\eps}}\) be \(L^2(0,1)\)-moderate regularizing \(f_1\) and \(f_2\). Let \(\net{b_\eps}\) be \(C([0,T], L^\infty(0,1))\)-moderate of log-type, \(\net{g_\eps}\) be \(L^2([0,T],H^{-2}(0,1))\)-moderate and let \(\net{c_\eps}\) be \(L^\infty(0,1)\)-moderate regularizing $b,g$ and $c$, respectively.
        Assume additionally that there exists a constant $c_0>0$ such that for every sufficiently small \(\eps>0\), it holds that 
		\begin{align}
		0<c_0\leq c_\eps(x), \quad \text{for all } x\in [0,1].\label{eqPropVWSExists}
		\end{align}
		Then there exists a very weak solution to Problem \ref{problemEB}.
	\end{thm}
	\pr
	By the respective moderateness of the regularizing nets, there are positive constants
    \(C_{f_1},C_{f_2}, C_c, C_b, C_g, N_{f_1}, N_{f_2}, N_c,N_g\) such that
    \begin{gather*}
    \|f_{1,\eps}\|_{H^2(0,1)} \leq C_{f_1} \eps ^{-N_{f_1}},\\
    \|f_{2,\eps}\|_{L^2(0,1)} \leq C_{f_2} \eps ^{-N_{f_2}},\\
    \|c_\eps\|_{L^\infty(0,1)} \leq C_c \eps^{-N_c}, \\
    \|b_\eps\|_{L^\infty(X_T)}, \leq C_b\log\dfrac{1}{\eps},\\
    \|g_\eps\|_{L^2([0,T], H^{-2}(0,1))} \leq C_g \eps^{-N_g}.
    \end{gather*}
For fixed \(\eps >0\), all conditions of Theorem \ref{WeakSolution} are verbatim fulfilled. The \textit{net of regularized problems} corresponding to Problem 1 is given by
\begin{align}
\partial_x^2\left(c_{\varepsilon}(x) \partial_x^2 u_{\varepsilon}(x,t) \right) + b_{\varepsilon}(x,t) \partial_x^2 u_{\varepsilon}(x,t) +\partial_t^2 u_{\varepsilon}(x,t) = g_{\varepsilon}(x,t), \quad (x,t) \in X_T,\\
u_{\varepsilon}(x,0) = f_{1,_{\varepsilon}}(x), \quad \partial_t u_{\varepsilon}(x,0) = f_{2,_{\varepsilon}}(x), \quad x \in [0,1], \\
 u_{\varepsilon}(0,t) = \partial_x u_{\varepsilon}(0,t) = 0,\quad u_{\varepsilon}(1,t) = \partial_x u_{\varepsilon}(1,t) = 0.
\end{align}
Thus, for every $\varepsilon>0$, we find a unique weak solution \(u_\eps\) to Problem 1. We claim that the net of solutions \(\net{u_\eps}\) is \(W\)-moderate. We use the energy estimate \eqref{EnergyEstimate2} of Theorem \ref{WeakSolution}.
It holds that
\begin{align*}
\|u_\eps\|^2_W &\leq \left( D^\eps_T \|f_{1,\eps}\|^2_{H^2(0,1)} + \|f_{2,\eps}\|^2_{L^2(0,1)} + \|g_\eps\|^2_{L^2([0,T],H^{-2}(0,1))}	\right) \exp(T F^\eps_T).
\end{align*}
Further, there exist constants \(C_D, C_F, N_F\geq 0\) such that
\begin{align}\begin{split}\label{eqDTEstimate}
D^\eps_T  &= \left(\|c_\eps\|_{L^\infty(0,1)} + c_0C_{1/2} (1+T)\right)/\min\left(\dfrac{c_0}{2},1\right)\\
&\leq C_c \eps^{-N_c}\left(1+c_0C_{1/2}(1+T)\right)/\min\left(\dfrac{c_0}{2},1\right) \\
&= C_D \eps^{-N_c},
\end{split}
\end{align}
and
\begin{align}
\begin{split}\label{eqFTEstimate}
\exp\left(T F^\eps_T\right) &= \exp\left(T \left(\|b_\eps\|_{L^\infty(X_T)} + 1 + c_0C_{1/2}(1+T)\right)/\min\left(\dfrac{c_0}{2},1\right)\right)\\
&\leq \exp\left(\dfrac{T\left(1 + c_0 C_{1/2}(1+T)\right)}{{\min\left(\dfrac{c_0}{2},1\right)}}\right) \exp\left(\dfrac{T C_b \log \dfrac{1}{\eps}}{\min\left(\dfrac{c_0}{2},1\right)}\right)\\
&\leq \exp\left(\dfrac{T\left(1 + c_0 C_{1/2}(1+T)\right)}{{\min\left(\dfrac{c_0}{2},1\right)}}\right) \eps^{\dfrac{-T C_b }{\min\left(\dfrac{c_0}{2},1\right)}}
\\
&=: C_F \eps^{-N_F}.
\end{split}
\end{align}
Finally, we obtain the following estimate
\begin{align}
\|u_\eps\|^2_W &\leq  \left( D^\eps_T \|f_{1,\eps}\|^2_{H^2(0,1)} + \|f_{2,\eps}\|^2_{L^2(0,1)} + \|g_\eps\|^2_{L^2([0,T],H^{-2}(0,1))}	\right) \exp\left(T F^\eps_T\right)\label{estimate} \\
&\leq \left(C_D C^2_{f_1}\eps^{-N_c-2N_{f_1}} + C^2_{f_2} \eps^{-2N_{f_1}} + C_g^2 \eps^{-2 N_g}\right) C_F \eps^{-N_F},\nonumber\\
&\leq C \eps^{-N},\nonumber
\end{align}
for some constants \(C, N>0\).
This proves that \(\net{u_\eps}\) is \(W\)-moderate.\ep

\begin{remark}\label{remModerateSolution}
By tracking the constants involved in the proof of Theorem \ref{veryWeakSolutionExists} we obtain the exact moderateness estimate for \(\net{u_\eps}\)
\begin{equation*}
    \|u_\eps\|_{W} \leq C_u \eps^{-N_u}, \label{eqModerateSolution2}
\end{equation*}
with constants \(C_u, N_u>0\) given by
\begin{gather}
C_u^2 = C_g^2 + C^2_{f_2} + C_c C^2_{f_1} (1+c_0C_{1/2}(1+T))/\min\left(\dfrac{c_0}{2}, 1\right),\label{eqExactConstant1}\\ 
2N_u = \max(N_c+2N_{f_1}, 2N_{f_2} , 2N_g) + C_b T / \min\left(\frac{c_0}{2},1\right).\label{eqExactConstant2}
\end{gather}
\end{remark}

\begin{remark}
It is important to note that Theorem \ref{veryWeakSolutionExists} does not cover all cases for which data regularizations correspond to a very weak solution. For example, assume that \(\net{c_\eps}\) and \(\net{f_{1,\eps}}\) are such that \(D^\eps_T \|f_{1,\eps}\|^2_{H^2(0,1)}\) is moderate, then the error estimate \eqref{estimate} still yields the moderateness of the sequential solution. This includes cases where \(\|c_\eps\|_{L^\infty(0,1)}\) has faster than moderate growth rate, while  \(\|f_{1,\eps}\|^2_{H^2(0,1)}\) is a zero sequence that vanishes quickly, and vice-versa. There are many other cases including zero sequences. However, a general study of these cases is of less importance since they highly depend on the specific growth and vanishing rates of the data regularizations. 

It is most natural to consider the case where none of the distributional data are zero, and therefore none of the regularizations and their estimates are zero sequences.
In this case, Theorem \ref{veryWeakSolutionExists} provides the existence results of a very weak solution while considering optimally large sets of regularizations for each of the data. Thus, let us further assume that none of the regularization estimates are zero sequences.
\end{remark}
 
\begin{remark}
Theorem \ref{veryWeakSolutionExists} can be improved by considering a sequence of positive lower bounds \(c_{0,\eps}\) of \(c_\eps\), i.e.,
\begin{equation*}
0<c_{0,\eps}\leq c_{\varepsilon}(x), \quad x \in (0,1),
\end{equation*} as it is not necessary that $c_0$ is independent of epsilon. Assume that we are not considering the case when the first factor on the RHS of Equation \eqref{EnergyEstimate2} is a zero sequence. Then, in order to obtain \(W\)-moderateness of \(\net{u_\eps}\), we need a moderate estimate of the term \(\exp(TF^\eps_T)\), as in \eqref{eqFTEstimate}. 
This will be the case whenever \(\dfrac{1}{c_{0,\eps}}\), \(||b_\eps||_{L^\infty(X_T)}\) and \(\dfrac{||b_\eps||_{L^\infty(X_T)}}{c_{0,\eps}}\) are moderate of log-type. Let us consider the following two extreme cases of the previous conditions. Firstly, if \(c_{0,\eps} = c_0\), as in Theorem \ref{veryWeakSolutionExists}, then \(\net{b_\eps}\) can be any \(C([0,T], L^\infty(0,1))\) moderate net of log-type. This case will be considered in the rest of the article. Secondly, let \(\net{b_\eps}\) be \(C([0,T], L^\infty(0,1))\)-bounded. Then it can be easily seen that the results hold  for any sequence of positive lower bounds \(\net{c_{0,\eps}}\) of \(L^\infty(0,1)\)-moderate net \(\net{c_\eps}\), satisfying
\begin{equation*}
0<\dfrac{M}{\log(1/\eps)}\leq c_{0,\eps} \leq c_\eps(x), \quad x \in (0,1),
\end{equation*}
for some \(M>0\). 
\end{remark}
	
	Now we give an example of a possible construction of regularizations of given distributional coefficients and data \eqref{eqDistData} which satisfy assumptions of Theorem \ref{veryWeakSolutionExists}. First, we recall a few definitions and results \cite{T1967} that will play an important role in the sequel. For the completeness of this work, we present a slightly modified proof of \cite[Thm. 21.1]{T1967} in Lemma \ref{lem3.5}.
	
\begin{definition}~
\begin{itemize}
\item[(a)] Let $C_c(\Omega)$ be the space of continuous functions on an open subset $\Omega\subseteq \mathbb{R}^d$ with compact support. A \textbf{Radon measure} is a linear functional on \(C_c(\Omega)\) such that for any compact subset \(K\subset \Omega\), there is a constant \(M_K\geq 0\) such that, for all continuous functions $\phi$ vanishing identically outside of $K$, there holds
	\[
	|\langle f, \phi \rangle | \leq M_K \sup \limits_{x \in K}|\phi(x)|.
	\]
\item[(b)] Distribution \(f \in \D'(\R^d)\) is a \textbf{positive distribution} 
if it satisfies
	\[
	\langle f, \varphi \rangle  \geq 0, 
	\]
	for every \(\varphi \in \D(\R^d)\) such that \(\varphi(x) \geq 0\).
\end{itemize}
\end{definition}

	\begin{lemma}
	\label{lem3.5}
	    A positive distribution is a positive Radon measure.
	\end{lemma}
    \pr
        Consider \(f \in \D'(\Omega)\) and let \(K\subseteq \Omega\) an arbitrary compact set.
        Choose \(\eta\in \D(\Omega)\) such that \(\eta(x) = 1\) on \(K\). For any real-valued \(\phi\in \D(K)\) we have
		\[
		-\|\phi\|_{L^\infty(0,1)} \eta(x) \leq \phi(x) \leq \|\phi\|_{L^\infty(0,1)} \eta(x), \quad \forall x \in \Omega.
		\]
		Since the the functions \(\phi + \|\phi\|_{L^\infty(0,1)} \eta \) and \(\|\phi\|_{L^\infty(0,1)} \eta - \phi\) are positive we have
		\begin{align*}
		    \langle f,\phi + \|\phi\|_{L^\infty(0,1)} \eta \rangle \geq 0,\quad \langle f, \|\phi\|_{L^\infty(0,1)} \eta - \phi\rangle \geq 0.
		\end{align*}
		We obtain 
		\[
		-\|\phi\|_{L^\infty(0,1)} \langle f, \eta\rangle \leq \langle f,\phi\rangle \leq \|\phi\|_{L^\infty(0,1)}\langle f, \eta\rangle,
		\]
        and therefore
		\[
		 |\langle f, \phi\rangle| \leq  \|\phi\|_{L^\infty(0,1)}\langle f, \eta\rangle,
		\]
        which means that \(f\) is a distribution of order zero. By \cite[Proposition 21.2]{T1967} \(f\) is a Radon measure.    \ep

\begin{lemma}
\label{lemdist}
	    Let \(c \in \D'(0,1)\) and \(c_0>0\). Then there exists a regularization \(\net{c_\eps}\) of \(c\) such that 
\begin{equation}
\label{31dist}
	    0<c_0\leq c_\varepsilon(x), \quad \forall x \in [0,1],
\end{equation}	   
	    if and only if there exists a positive distribution \(\mu\) such that \(c = c_0 + \mu\).
	\end{lemma}
	\pr
Let \(c\in \D'(0,1)\), \(c_0>0\) and \(c_\eps\) such that \(0<c_0\leq c_\varepsilon(x)\).
        We have that
		\[
		c_0 \int_0^1 \varphi(x)dx \leq \langle c_\eps,\varphi\rangle,
		\]
		if \(\varphi\in \D(0,1)\) with \(\varphi \geq 0\) and by passing to the limit
		\[
		c_0 \int_0^1 \varphi(x)dx \leq \langle c, \varphi \rangle.
		\]
		Set \(\mu = c - c_0\). Then it is clear that \(\mu\) is a positive distribution.\\
Let now \(c_0>0\) and \(c = c_0 + \mu\in \D'(0,1)\) with \(\mu\in \D'(0,1)\) a positive distribution. Consider the extension \(\tilde{c}\in \D'(\R)\) given by \(\tilde{c} = c_0 + \mu\) on the whole of \(\R\). Any regularization of \(\tilde{c}\) is a regularization of \(c\). Let \(\varphi\) be a Friedrichs mollifier and define \(c_\eps = \tilde{c} * \varphi_\eps\). We estimate
		\begin{equation*}
		c_\eps(x) = c_0*\varphi_\eps + \mu * \varphi_\eps \geq  c_0 * \varphi_\eps = c_0 \int_{\R}\varphi_\varepsilon(x+y) dy = c_0,
		\end{equation*}
		since \(\mu\) is a positive distribution, so \(\mu  *\varphi_\eps\) is non-negative. This concludes the proof.
	\ep

\begin{thm}[Regularizations]
	Let the initial values \(f_1, f_2\in \D'(0,1)\). Let the right-hand side \(g \in \D'(X_T)\). Let the coefficient \(b \in \D'(X_T)\) and let \(c \in \D'(0,1)\) of the form \(c = c_0 + \mu_c\), where \(c_0>0\) and \(\mu_c\in \D'(0,1)\) is a positive Radon measure. Then there exists a very weak solution to Problem \ref{problemEB}.
\end{thm}

\pr
We reduce the proof of the proposition, by constructing explicit regularizing nets of the data that satisfy the conditions in Theorem \ref{veryWeakSolutionExists}. Let \(\psi \in C^\infty_0(\R)\) be a Friedrichs mollifier. By Remark \ref{reModerateMollification}, the regularizing nets
		\begin{gather*}
		\net{f_{1,\eps}} = \net{f_1*\psi_\eps},\qquad \net{f_{2,\eps}} = \net{f_2*\psi_\eps},
		\end{gather*}
are \(H^2_0(0,1)\)-moderate and \(L^2(0,1)\)-moderate, respectively. Let $(c_{\varepsilon})_{\varepsilon}$ be any regularization of $c$ obtained by the construction explained in the proof of Lemma \ref{lemdist}. Then, we have the guarantee that the regularizing net $(c_{\varepsilon})_{\varepsilon}$ satisfies the uniform positivity condition \eqref{31dist} with the constant $c_0$. Let \(\net{\varphi_\eps}\subset C^\infty_0(\R^2)\) be a Friedrichs mollifier. We define the regularizing net
\begin{gather*}
g_\eps = g* \varphi_\eps|_{X_T},
\end{gather*}
which is \(L^2(X_T)\)-moderate by Remark \ref{reModerateMollification}, thus in particular \(L^2([0,T], H^{-2}(0,1))\)-moderate.
For \(b\) we use the regularizing net \(\net{b_\eps}\) defined by
\[
b_\eps = b * \varphi_{\lambda(\eps)}|_{X_T},
\]
with the reparametrisation function \(\lambda(\eps) = \left(\log \log\dfrac{1}{\eps} \right)^{-1}\). It follows that \(\net{b_\eps}\) is an \(L^\infty(X_T)\)-moderate net of log-type.
The regularizations \(\net{b_\eps}\), \(\net{g_\eps}\) and \(\net{c_\eps}\) satisfy the conditions of Theorem \ref{veryWeakSolutionExists}. Thus we find a \(W\)-moderate net \(\net{u_\eps}\) of solutions. We conclude that \(\net{u_\eps}\) is a very weak solution to \eqref{eqEB}.
\ep

By the above construction of regularizing nets we can guarantee the existence of a very weak solution. However this construction of regularizations by smooth nets is not applicable for the consistency Theorem \ref{thmConsistency}, in the case of non-continuous \(L^\infty\)-functions \(b\) and \(c\). To resolve this, we refer to Remark \ref{RemarkConsistencyRegularizations}.
	
	\section{Uniqueness and consistency results}\label{secConsistencyUniqueness}
	For each choice of regularizing nets which satisfy the conditions of Theorem \ref{veryWeakSolutionExists}, one obtains a unique sequential solution \(\net{u_\eps}\). However since the choice of regularizing nets is arbitrary, the resulting sequential solutions can still be different. 
		
\subsection{Uniqueness} As previously mentioned, we will formulate the uniqueness results in the neglibility sense. In other words, we seek conditions under which the coefficients and initial conditions are close enough, such that we can guarantee that the resulting sequential solutions are close enough. In the case of a very weak solution, we find stability under negligible deviations of the regularizing nets.
	
	\begin{thm}[Uniqueness]
		\label{thmVWSUniqueness}
		Let \(f_1\in \D'(0,1)\), \(f_2\in \D'(0,1)\), \(b \in \D'(X_T)\), \(g \in \D'(X_T)\) and \(c\in \D'(0,1)\). Suppose \(c\in \D'(0,1)\) is of the form \(c_0 +\mu\) for a positive constant \(c_0>0\) and for a positive Radon measure \(\mu\in \D'(0,1)\).
		Let \(\net{u_\eps}\) and \(\net{\tilde{u}_\eps}\) be two very weak solutions to the problem obtained through Theorem \ref{WeakSolution} with respective regularizing nets \(\net{f_{1,\eps}}\), \(\net{f_{2,\eps}}\), \(\net{b_\eps}, \net{c_\eps}, \net{g_\eps}\) and \(\net{\tilde{f}_{1,\eps}}\), \(\net{\tilde{f}_{2,\eps}}\)\(\net{\tilde{b}_\eps},\) \(\net{\tilde{c}_\eps}\), \(\net{\tilde{g}_\eps}\) which satisfy the following:
		
	\begin{itemize}
		\item \(\net{f_{1,\eps}}\) and \(\net{\tilde{f}_{1,\eps}}\) are regularizations of \(f_1\) that are \(H^2_0(0,1)\)-moderate, such that \(\net{f_{1,\eps} - \tilde{f}_{1,\eps}}\) is \(H^2(0,1)\)-negligible.
		
		\item \(\net{f_{2,\eps}}\) and \(\net{\tilde{f}_{2,\eps}}\) are regularizations of \(f_2\) that are \(L^2(0,1)\)-moderate, such that \(\net{f_{2,\eps} - \tilde{f}_{2,\eps}}\) is \(L^2(0,1)\)-negligible.
		
		\item \(\net{b_\eps}\) and \(\net{\tilde{b}_\eps}\) are regularizations of \(b\) that are \(C([0,T], L^\infty(0,1))\)-moderate of log-type, such that \(\net{b_\eps-\tilde{b}_\eps}\) is \(L^2([0,T],L^\infty(0,1))\)-negligible.
		
		\item \(\net{c_\eps}\) and \(\net{\tilde{c}_\eps}\) are regularizations of \(c\) that
		are \(L^\infty(0,1)\)-moderate, that satisfy the bounds
		\[
		0<c_0\leq c_\eps(x),\quad c_0\leq \tilde{c}_\eps(x), \quad \text{for all } x \in [0,1],
		\]
        and such that \(\net{c_\eps-\tilde{c}_\eps}\) is \(L^\infty(0,1)\)-negligible.
		
		\item \(\net{g_\eps}\) and \(\net{\tilde{g}_\eps}\) are regularizations of \(g\) that are \(L^2([0,T], H^{-2}(0,1))\)-moderate, such that \(\net{g_\eps - \tilde{g}_\eps}\) is \(L^2([0,T], H^{-2}(0,1))\)-negligible.
	\end{itemize}
	Then \(\net{u_\eps-\tilde{u}_\eps}\) is \(W\)-negligible.
	\end{thm}
	\pr
	Consider solution nets \(\net{u_\eps}\) and \(\net{\tilde{u}_\eps}\) and regularizing nets of the data \(\net{f_{1,\eps}}\), \(\net{\tilde{f}_{1,\eps}}\), \(\net{f_{2,\eps}}\), \(\net{\tilde{f}_{2,\eps}}\), \(\net{b_\eps}\), \(\net{\tilde{b}_\eps}\), \(\net{c_\eps}\), \(\net{\tilde{c}_\eps}\), \(\net{g_\eps}\) and \(\net{\tilde{g}_\eps}\) that satisfy the conditions of Theorem \ref{thmVWSUniqueness}.
	We need to prove that the net of differences \(\net{u_\eps-\tilde{u}_\eps}\) is \(W\)-negligible.
	The weak formulations of the regularized equations are
	\begin{align}
		\langle c_\eps(x) \partial_x^2 u_\eps(x,t),\partial_x^2 v(x)\rangle  + \langle b_\eps(x,t) \partial_x^2 u_\eps(x,t), v(x)\rangle + \langle \partial_t^2 u_\eps(x,t), v(x)\rangle  = \langle g_\eps(x,t),v(x)\rangle, \label{eqEBWeakReg1}\\
		\langle \tilde{c}_\eps(x) \partial_x^2 \tilde{u}_\eps(x,t),\partial_x^2 v(x)\rangle  + \langle\tilde{b}_\eps(x,t) \partial_x^2 \tilde{u}_\eps(x,t), v(x)\rangle + \langle \partial_t^2 \tilde{u}_\eps(x,t), v(x)\rangle  = \langle \tilde{g}_\eps(x,t),v(x)\rangle, \label{eqEBWeakReg2}
		\end{align}
		for all \( t \in [0,T]\) and \(v \in H^2_0(0,1)\).
		Subtracting equation \eqref{eqEBWeakReg2} from \eqref{eqEBWeakReg1} and rearranging the terms gives
		\begin{align}
		\begin{split}\label{eqUniquenessWeakDiff}
		&\langle c_\eps(x) \partial_x^2 (u_\eps-\tilde{u}_\eps), \partial_x^2 v(x) \rangle + \langle b_\eps(x,t) \partial_x^2 (u_\eps-\tilde{u}_\eps), v(x)\rangle + \langle \partial_t^2 (u_\eps-\tilde{u}_\eps), v(x)\rangle
		 \\
		&= \langle g_\eps(x,t)-\tilde{g}_\eps(x,t),v(x)\rangle
		+ \langle (\tilde{c}_\eps(x)-c_\eps(x))\partial_x^2 \tilde{u}_\eps, \partial_x^2 v(x)\rangle + \langle(\tilde{b}_\eps(x,t)-b_\eps(x,t)) \partial_x^2 \tilde{u}_\eps, v(x) \rangle.
		\end{split}
		\end{align}
		For each \(\eps \in (0,1]\) and \(t \in [0,T]\), the right-hand side of \eqref{eqUniquenessWeakDiff} defines a functional \(h_\eps(x,t) \in H^{-2}(0,1)\) by
		\[
		h_\eps(x,t) := g_\eps(x,t)-\tilde{g}_\eps(x,t) + \partial_x^2(\tilde{c}_\eps(x)-c_\eps(x))\partial_x^2 \tilde{u}_\eps(x,t) + (\tilde{b}_\eps(x,t)-b_\eps(x,t)) \partial_x^2 \tilde{u}_\eps(x,t).
		\]
		It can be seen that \(h_\eps\) is in \(L^2([0,T],H^{-2}(0,1))\) since \(g_\eps, \tilde{g}_\eps \in L^2([0,T], H^{-2}(0,1))\), \(c_\eps, \tilde{c}_\eps \in L^\infty(0,1)\), \(b_\eps, \tilde{b}_\eps \in L^\infty(X_T)\) and \(\tilde{u}_\eps\in L^\infty([0,T],H^2(0,1))\). From Theorem \ref{WeakSolution} it follows that \(u_\eps - \tilde{u}_\eps\) is the unique function in \(L^2([0,T], H^2_0(0,1))\) that satisfies \eqref{eqUniquenessWeakDiff} and the initial conditions 
		\[
		(u_\eps-\tilde{u}_\eps)(x,0) = f_{1,\eps} - \tilde{f}_{1,\eps}, \quad \text{and}\quad \partial_t (u_\eps-\tilde{u}_\eps)(x,0) = f_{2,\eps} - \tilde{f}_{2,\eps}.
		\]
		The energy estimate \eqref{EnergyEstimate2} still applies and leads to:
		\begin{align*}
		\left\|u_\eps-\tilde{u}_\eps\right\|^2_W &\leq \left(D^\eps_T \|f_{1,\eps} - \tilde{f}_{1,\eps}\|^2_{H^2(0,1)} + \|f_{2,\eps} - \tilde{f}_{2,\eps}\|^2_{L^2(0,1)} + \|h_\eps\|^2_{L^2([0,T],H^{-2}(0,1))}\right) \exp\left(T F^\eps_T\right).
		\end{align*}
		As in \eqref{eqDTEstimate} and \eqref{eqFTEstimate} the constants \(C_D\) and \(C_F, N_F\) satisfy
		\begin{gather*}
		    D_T^\eps \leq C_D\,\varepsilon^{-N_c},\qquad \text{ and }\qquad
		    \exp\left(T F_T^\eps\right) \leq C_F\, \eps^{-N_F}.
		\end{gather*}
	    On the other hand, for any \(q>0\), from the negligibility condition, there exist positive constants  \(C_{q}^{f_1}\) and \(C_{q}^{f_2}\) such that as $\varepsilon\to 0$, there holds
		\begin{gather}
		\|f_{1,\eps}-\tilde{f}_{1,\eps}\|_{H^2(0,1)} \leq C^{f_1}_{q} \eps^{q},\\
		\|f_{2,\eps}-\tilde{f}_{2,\eps}\|_{L^2(0,1)} \leq C^{f_2}_{q} \eps^{q},
		\end{gather}
		and as we will verify shortly there is \(C^h_{q}>0\) such that
\begin{equation}
		\|h_\eps\|_{L^2([0,T],H^{-2}(0,1))} \leq C^h_{q} \eps^{q}, \quad \text{as }\eps \to 0. \label{eqhepsneg}
\end{equation}
	    First, for any \(q>0\) we obtain 
		\begin{align*}
		\left\|u_\eps-\tilde{u}_\eps\right\|^2_W &\leq \left(D^\eps_T \|f_{1,\eps} - \tilde{f}_{1,\eps}\|^2_{H^2(0,1)} + \|f_{2,\eps} - \tilde{f}_{2,\eps}\|^2_{L^2(0,1)} + \|h_\eps\|^2_{L^2([0,T],H^{-2}(0,1))}\right) \exp\left(T F^\eps_T\right)\\
		&\leq \left(C_D \eps^{-N_c} \left(C_{\frac{q+N_c}{2}}^{f_1} \eps^{\tfrac{q+N_c}{2}}\right)^2 + \left(C_{q/2}^{f_2} \eps^{\tfrac{q}{2}}\right)^2 + \left(C_{q/2}^h \eps^{\tfrac{q}{2}}\right)^2 \right) C_F \eps^{-N_F}\\
		&= C_F\left(C_D \left(C_{\frac{q+N_c}{2}}^{f_1}\right)^2 + \left(C_{q/2}^{f_2}\right)^2 + \left(C_{q/2}^h\right)^2\right) \eps^{q - N_F}\\
		&= C^2_{\tfrac{q-N_F}{2}} \eps^{q-N_F},
		\end{align*}
and thus we have that \(\net{u_\eps - \tilde{u}_\eps}\) is \(W\)-negligible.

Finally, let us verify the estimate \eqref{eqhepsneg}, i.e., let us show that \(h_\eps\) is \(L^2([0,T],H^{-2}(0,1))\)-negligible. Let \(q>0\) arbitrary. 		First, we can easily see that
		\begin{align*}
		\|h_\eps\|_{L^2([0,T],H^{-2}(0,1))} &\leq \|g_\eps - \tilde{g}_\eps\|_{L^2([0,T],H^{-2}(0,1))} + \|\partial_x^2\left((\tilde{c}_\eps-c_\eps)\partial_x^2 \tilde{u}_\eps\right) \|_{L^2([0,T],H^{-2}(0,1))} \\
		&\quad+ \| (\tilde{b}_\eps - b_\eps) \partial_x^2 \tilde{u}_\eps\|_{L^2([0,T],H^{-2}(0,1))}
		\end{align*}
It is thus sufficient that each of the terms is \(L^2([0,T],H^{-2}(0,1))\)-negligible.
    	By assumption, there is \(C_q^g>0\) such that
		\[
		\|g_\eps-\tilde{g}_\eps\|_{L^2([0,T],H^{-2}(0,1))} \leq C_q^g \eps^q,\quad \text{as }\eps \to 0.
		\] 
		Second, from H\"older's inequality it follows that 
		\begin{align*}
		\|(\tilde{b}_\eps - b_\eps)\partial_x^2\tilde{u}_\eps\|_{L^2([0,T], H^{-2}(0,1))} &\leq \|(\tilde{b}_\eps - b_\eps)\partial_x^2\tilde{u}_\eps\|_{L^2([0,T], L^2(0,1))}\\
		&\leq \|(\tilde{b}_\eps - b_\eps)\|_{L^2([0,T],L^\infty(0,1))} \|\tilde{u}_\eps\|_{L^\infty([0,T],H^2(0,1))} \\
		&\leq C_{q+N_{\tilde{u}}}^b \eps^{q+ N_{\tilde{u}}}  C_{\tilde{u}} \eps^{-N_{\tilde{u}}}\\
		&= C_{q+N_{\tilde{u}}}^b C_{\tilde{u}}\eps^q,
		\end{align*}
		where \(C_{\tilde{u}}, N_{\tilde{u}}\) satisfy \eqref{eqExactConstant1} and \eqref{eqExactConstant2} and 
		\[
		 \|(\tilde{b}_\eps - b_\eps)\|_{L^2([0,T],L^\infty(0,1))} \leq C_q^b \eps^q,\quad \text{as }\eps \to 0.
		\]
		Lastly, if
		\[
		 \|\tilde{c}_\eps-c_\eps\|_{L^\infty(0,1)} \leq C^c_q \eps^q,\quad \text{as }\eps \to 0,
		\]
		then the Leibniz rule and H\"older inequality give
		\begin{align*}
		\|\partial_x^2(\tilde{c}_\eps-c_\eps)\partial_x^2\tilde{u}_\eps\|_{L^2([0,T],H^{-2}(0,1))} &\leq 
		\|(\tilde{c}_\eps-c_\eps)\partial_x^2\tilde{u}_\eps\|_{L^2(X_T)}\\
		&\leq \|\tilde{c}_\eps-c_\eps\|_{L^\infty(0,1)}\|\tilde{u}_\eps\|_{L^\infty([0,T],H^2(0,1))}\\
		&\leq C_{q+N_{\tilde{u}}}^c \eps^{q + N_{\tilde{u}}} C_{\tilde{u}}\eps^{-N_{\tilde{u}}}\\
		&= C_{q+N_{\tilde{u}}}^c C_{\tilde{u}}\eps^{q},
		\end{align*} which concludes the proof.
	\ep
	
	\begin{remark}
	One can consider a different but similar notion of uniqueness than Theorem \ref{thmVWSUniqueness}. In this notion we aim to conclude \(\|u_\eps - \tilde{u}_\eps\|_W\to 0\). With the same notation as above, assume that moderateness-estimates of the data is given in terms of the same constants appearing in the proof of the theorem. If we assume
	\begin{align}
	    \|f_{1,\eps}-\tilde{f}_{1,\eps}\|_{H^2(0,1)} &= o\left( \eps^{\tfrac{N_c + N_F}{2}}\right)\label{eqNumUniquef1}
	    \\
	    \|f_{2,\eps}-\tilde{f}_{2,\eps}\|_{L^2(0,1)} &= o\left( \eps^{\tfrac{N_F}{2}}\right),\label{eqNumUniquef2}\\
	    \|g_\eps -\tilde{g}_\eps\|_{L^2([0,T], H^{-2}(0,1))} &= o\left(\eps^{\tfrac{N_F}{2}}\right),\\
        \|b_\eps - \tilde{b}_\eps\|_{L^2([0,T], L^\infty(0,1))} &= o\left(\eps^{N_{\tilde{u}}+N_F/2}\right),\\
        \|c_\eps -\tilde{c}_\eps\|_{L^\infty(0,1)} &= o\left( \eps^{N_{\tilde{u}} + N_F/2}\right),\label{eqNumUniquec}
	\end{align}
    then a precise inspection of the proof above immediately shows that 
   \begin{align*}
		\left\|u_\eps-\tilde{u}_\eps\right\|^2_W &\leq \left(D^\eps_T \|f_{1,\eps} - \tilde{f}_{1,\eps}\|^2_{H^2(0,1)} + \|f_{2,\eps} - \tilde{f}_{2,\eps}\|^2_{L^2(0,1)} + \|h_\eps\|^2_{L^2([0,T],H^{-2}(0,1))}\right) \exp\left(T F^\eps_T\right)\\
		&\leq \left(C_D \eps^{-N_c} o\left(\eps^{\tfrac{N_c +N_F}{2}}\right)^2 + o\left(\eps^{\tfrac{N_F}{2}}\right)^2 + o\left(\eps^{\tfrac{N_F}{2}}\right)^2 \right) C_F \eps^{-N_F}\\
		&= o(1).
		\end{align*}
	This observation is useful for numerical methods, where it is practically difficult to consider negligible differences in regularizing nets. 
	\end{remark}

	\subsection{Consistency}
	The very weak solution obtained in Theorem \ref{veryWeakSolutionExists} coincides with the weak solution guaranteed by Theorem \ref{WeakSolution}, when the data are sufficiently regular functions. This results will be more formally treated in the following theorem.
	
	\begin{thm}[Consistency]
	\label{thmConsistency}
	    Let \(f_1 \in H^2_0(0,1)\) and \(f_2 \in L^2(0,1)\).
		Let \(b \in C([0,T],L^\infty(0,1))\), \(g\in L^2([0,T],H^{-2}(0,1))\) and \(c \in L^\infty(0,1)\) with
		\begin{gather*}
		0<c_0\leq c(x), \quad c_0 \in \R_+,
		\end{gather*}
		for all \(x\in [0,1]\) and small \(\eps>0\).
		Let \(u\in L^2([0,T], H^2_0(0,1))\) be a weak solution to problem  \eqref{eqEB} and let \(\net{u_\eps}\subset L^2([0,T],H^2_0(0,1))\) be a very weak solution of \eqref{eqEB} with corresponding regularizing nets \(\net{f_{1,\eps}}\), \(\net{f_{2,\eps}}\), \(\net{b_\eps}, \net{c_\eps}\) and \(\net{g_\eps}\) where \( 0<c_0\leq c_\eps(x)\). Suppose that \(f_{1,\eps}\to f_1\) in \(H^2(0,1)\), \(f_{2,\eps}\to f_2\) in \(L^2(0,1)\), \(b_\eps \to b\) in \(L^2([0,T], L^\infty(0,1))\), \(g_\eps \to g\) in \(L^2([0,T], H^{-2}(0,1))\) and \(c_\eps \to c\) in \(L^\infty(0,1)\) as \(\eps \to 0\). 	Then \(\net{u_\eps}\) converges strongly to \(u\) in \(W\).
	\end{thm}
	\pr
		The weak solution \(u\) and the sequential solution \(\net{u_\eps}\) satisfy the following weak formulations respectively:
		\begin{gather*}
		\langle c(x) \partial_x^2 u(x,t),\partial_x^2 v(x)\rangle  + \langle b(x,t) \partial_x^2 u(x,t), v(x)\rangle + \langle \partial_t^2 u(x,t), v(x)\rangle  = \langle g(x,t),v(x)\rangle,,\\
		\langle c_\eps(x) \partial_x^2 u_\eps(x,t),\partial_x^2 v(x)\rangle  + \langle b_\eps(x,t) \partial_x^2 u_\eps(x,t), v(x)\rangle + \langle \partial_t^2 u_\eps(x,t), v(x)\rangle  = \langle g_\eps(x,t),v(x)\rangle, 
		\end{gather*} for all \( t \in [0,T]\) and \(v \in H^2_0(0,1)\), with the initial conditions
		\begin{align*}
		    u(x,0) = f_1(x),\quad \partial_t u(x,0) = f_2(x),\\
		    u_\eps(x,0) = f_{1,\eps}(x),\quad \partial_t u_\eps(x,0) = f_{2,\eps}(x).
		\end{align*}
		Then the difference of solutions \(u_\eps -u\) satisfies
		\begin{align}
		\begin{split}\label{eqConsistencyDiff}
		\langle c_\eps(x) \partial_x^2 (u_\eps-u), \partial_x^2 v(x) \rangle + \langle b_\eps(x,t) \partial_x^2 (u_\eps -u), v(x)\rangle + \langle \partial_t^2 (u_\eps -u), v(x)\rangle\\
		= \langle g_\eps(x,t)-g(x,t),v(x)\rangle
		+ \langle (c(x)-c_\eps(x))\partial_x^2 u, \partial_x^2 v(x)\rangle + \langle(b(x,t)-b_\eps(x,t)) \partial_x^2 u, v(x) \rangle,
		\end{split}
		\end{align}
		for \(t \in [0,T]\), \(v \in H^2_0(0,1)\), and the initial conditions 
		\begin{equation}
		(u_\eps - u)(x,0) = f_{1,\eps}(x)-f_1(x),\quad \partial_t(u_\eps-u)(x,0)= f_{2,\eps}(x) -f_2(x).
		\label{eqConsistencyIc}
		\end{equation}
		For each \(t \in [0,T]\) the right hand side of \eqref{eqConsistencyDiff} induces a functional \(h_\eps(x,t)\in H^{-2}(0,1)\). Moreover \(h_\eps\in L^2([0,T], H^{-2}(0,1))\) since \(g, g_\eps \in L^2([0,T], H^{-2}(0,1))\), \(c, c_\eps \in L^\infty(0,1)\), \(b, b_\eps \in L^\infty(X_T)\) and \(u\in L^\infty([0,T],H^2(0,1))\).
		
		By Theorem \ref{WeakSolution}, \(u_\eps - u\) is the unique solution to the weak formulation
		\begin{align*}
		\langle c_\eps(x) \partial_x^2 (u_\eps -u), \partial_x^2 v(x) \rangle + \langle b_\eps(x,t) \partial_x^2 (u_\eps - u), v(x)\rangle + \langle \partial_t^2 (u_\eps - u), v(x)\rangle
		= \langle h_\eps(x,t),v(x)\rangle,
		\end{align*}
		with the initial conditions \eqref{eqConsistencyIc}. The energy estimate \eqref{EnergyEstimate2} then applies. 
	    Notice that \(D_T^\eps\) and \(F_T^\eps\) are bounded as \(\eps \to 0\). We then have
		\begin{align*}
		\|u_\eps-u\|_W &\leq  \left(D_T^\eps \|f_{1,\eps}-f_1\|^2_{H^2(0,1)} + \|f_{2,\eps} - f_2\|^2_{L^2(0,1)} +  \|h_\eps\|^2_{L^2([0,T],H^{-2}(0,1))}\right)\exp{\left(T F_T^\eps\right)}  \\ 
		& = \left(D_T^\eps o(1) + o(1) +  o(1)\right)\exp{\left(T F_T^\eps\right)} \\
		&= o(1),
		\end{align*}
		as we can estimate $h_{\varepsilon}$ as follows:
		\begin{align*}
		\|h_\eps\|_{L^2([0,T],H^{-2}(0,1))} &\leq \|g_\eps - g\|_{L^2([0,T],H^{-2}(0,1))} + \|\partial_x^2\left((c-c_\eps)\partial_x^2 u\right) \|_{L^2([0,T],H^{-2}(0,1))} \\
		&\quad+ \| (b - b_\eps) \partial_x^2 u\|_{L^2([0,T],H^{-2}(0,1))}\\
		&\leq o(1) + \|c-c_\eps\|_{L^\infty(0,1)} \|u\|_{L^2(X_T)} + \|b-b_\eps\|_{L^2([0,T],L^\infty(0,1))} \|u\|_{L^\infty([0,T], H^2(0,1))}\\
		&= o(1) + o(1) O(1) + o(1)O(1)\\
		&= o(1).
		\end{align*}
		This concludes the proof that \(\|u_\eps-u\|_{W} \to 0\) as $\varepsilon\to 0$.
\ep
\begin{remark}~
\label{RemarkConsistencyRegularizations}
\begin{itemize}
\item[(a)] If \(b\) and \(c\) are not continuous, the regularizing nets \(\net{b_\eps}\) and \(\net{c_\eps}\) such that \(\|b_\eps - b\|_{L^\infty(X_T)}\to 0\) and \(\|c_\eps - c\|_{L^\infty(0,1)}\to 0\) cannot be constructed via mollification, as previously mentioned. This can be seen as follows. The regularizing functions \(c_\eps = c\,*\varphi_\eps\) are continuous. The property \(\|c_\eps - c\|_{L^\infty(0,1)}\to 0\) implies also convergence in \(C([0,1])\), where the essential supremum becomes just the ordinary supremum in case of continuous functions. Thus, \(c\) is in \(C([0,1])\), which contradicts the assumption that $c$ is not continuous. Similarly for \(b\).
\item[(b)] This problem can be handled in different ways. One possibility is to consider the regularizing nets \(b_\eps\) and \(c_\eps\) of \(L^\infty\) functions, abandoning constructions of regularizing nets via mollifiers and smooth functions. For example, for a given zero sequence \(\net{p_\eps}\to 0\) in \(L^\infty(0,1)\), one can consider the net \(\net{c_\eps}\) of the form \(\net{c+p_\eps}\). Similarly for \(b\).

Another possibility is to restrict the consistency theorem to continuous \(b,c\). Then \(L^\infty\)-convergence of regularizing nets can be obtained by the following construction.
Let \(c \in C([0,1])\) and let \(\tilde{c}\in C(\R)\) such that \(\tilde{c}|_{[0,1]} = c\). If \(\varphi\) is a Friedrichs mollifier and $c_\eps = \tilde{c} * \varphi_\eps$, then for any \(x \in [0,1]\),
\begin{align*}
\left|c_\eps(x)-c(x)\right| & = \left|\int_\R \left(\tilde{c}(x-y) - c(x)\right) \varphi_\eps(y) dy\right|\\
&\leq \int_\R \left|\tilde{c}(x-y)-c(x) \right|\varphi_\eps(y) dy\\
&\leq \sup\limits_{y \in [-\eps, \eps]} \left|\tilde{c}(x-y) - \tilde{c}(x) \right|\leq \epsilon
\end{align*}
whenever \(\eps\leq \delta_\epsilon\), where \(\delta_\epsilon\) is such that
\(|x-z|\leq \delta_\epsilon \implies |\tilde{c}(x) - \tilde{c}(z)| \leq \epsilon\) for any \(x,z \in [-1,2]\). 
Such \(\delta_\epsilon\) exists since \(\tilde{c}\) is continuous on \([-1,2]\) and therefore uniformly continuous. Note also that  interval $[-1,2]$ can be replaced with an arbitrary closed neighbourhood $K$ of $[0,1]$.
\end{itemize}
\end{remark}
	
\section{Numerical simulations}

In this section we aim to illustrate our theoretical results via numerical simulations \cite{quart}.  More precisely, we will test performance of the numerical schemes for solving Problem 1.

We describe briefly the method we use our simulations. Firstly, we use a hybrid method in which we discretize the space dimension and apply a finite element method, while in the time dimension we work with a finite difference method. We discretize the space and time dimensions in \(n\) and \(m\) parts, respectively with discretization steps \(h = 1/n\) and \(h_t = 1/m\). For the space discretization, we consider the finite element space \(V_h\subset H^2_0(0,1)\) of continuous functions \(f\) on \([0,1]\) that are cubic on intervals of the type \([ih,(i+1)h]\), \(i = 0,\dots,n-1\) and such that the boundary conditions \eqref{eqbc} hold. Note that without the boundary conditions, this space is known as \(V^{1,3}\). As our boundary conditions are equal to zero, we denote this space with $V^{1,3}_0$.

Let \(u\) denote a numerical solution and let \(u_{ex}\) be the exact solution, then in the following experiments we consider the following error:
\begin{equation}
E_{L^2}(u, u_{ex}) = \max\left\{\|u(x,t)-u_{ex}(x,t)\|_{L^2(0,1)}:\, t = ih_t, \,i = 0,\dots,m\right\}.
\label{errorL2}
\end{equation}
%

\subsection{Regularization }
As we described previously in Section \ref{sec:regular}, regularization plays the central role in our work. Now, we explain the numerical procedure for the regularization of the coefficients. The standard mollifier in \(\Rn\) is given by
\begin{equation}
\varphi(x)  = \begin{cases}
A \exp\left(\dfrac{1}{1-\|x\|^2}\right), &\text{if } \|x\|<1,\\
0, &\text{otherwise},
\end{cases}\label{eq: standard mollifier}
\end{equation}
where \(A\) is obtained such that
\[
\int_{\R^n} \varphi(x) dx = 1.
\]
Then \(\varphi_\varepsilon(x) = \dfrac{1}{\varepsilon^n} \varphi\left(\dfrac{x}{\varepsilon}\right)\) is a mollifying net.
In the one-dimensional case we will employ the standard mollifier as in \eqref{eq: standard mollifier}. For the ease of numerical simulations and analysis, instead of considering two-dimensional standard molifier, we employ the separable molifier of the form \(\rho(x,t) = \varphi(x)\varphi(t)\). This mollifier has the nice property that \(\delta(x)*\rho(x,t) = \varphi(t)\) and similarly \(\delta(t)*\rho(x,t) = \varphi(x)\). 

Additionally we use an extension procedure as in Lemma \ref{lemdist}
to mollify extensions of the distributions on their domain. This will mostly be restricted to a constant mollification of a constant function.

\subsection{Regular case}
First, as the consistency theorem guarantees, when the data are sufficiently nice functions, the very weak solution obtained by the existence theorem, coincides with the classical or weak solution. Thus, we will first consider a problem where the solution is analytically known. We will use this simple problem to investigate the convergence of the numerical methods we will apply later on irregular cases, i.e., when our data are irregular functions. 
\begin{figure}[H]
  \centering
\subfigure[]{
  \includegraphics[height=4.5cm]{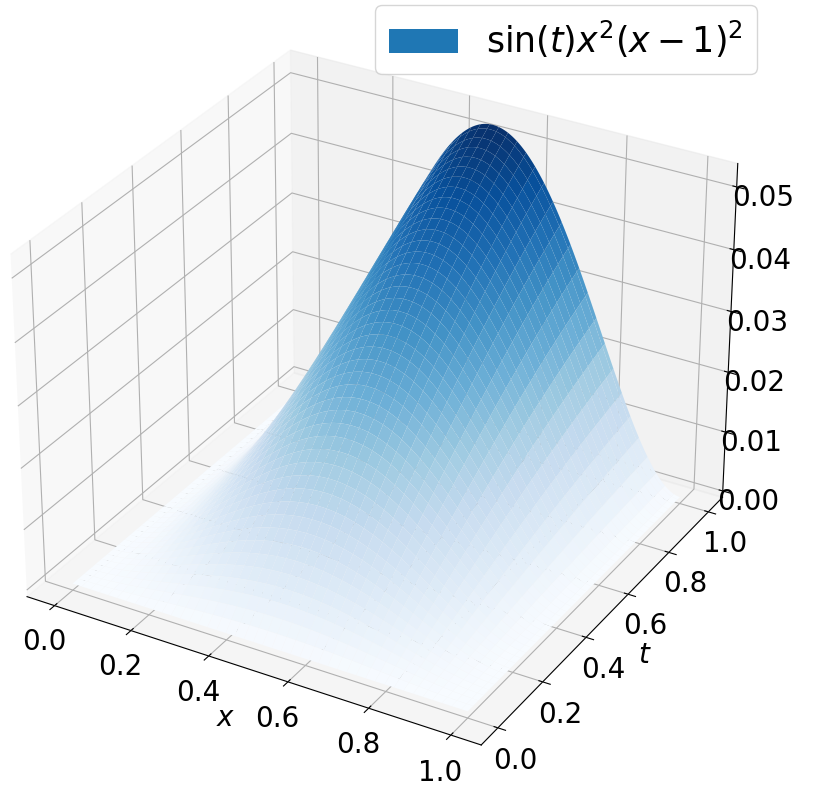}
  }
  \subfigure[]{
    \includegraphics[height=4.5cm]{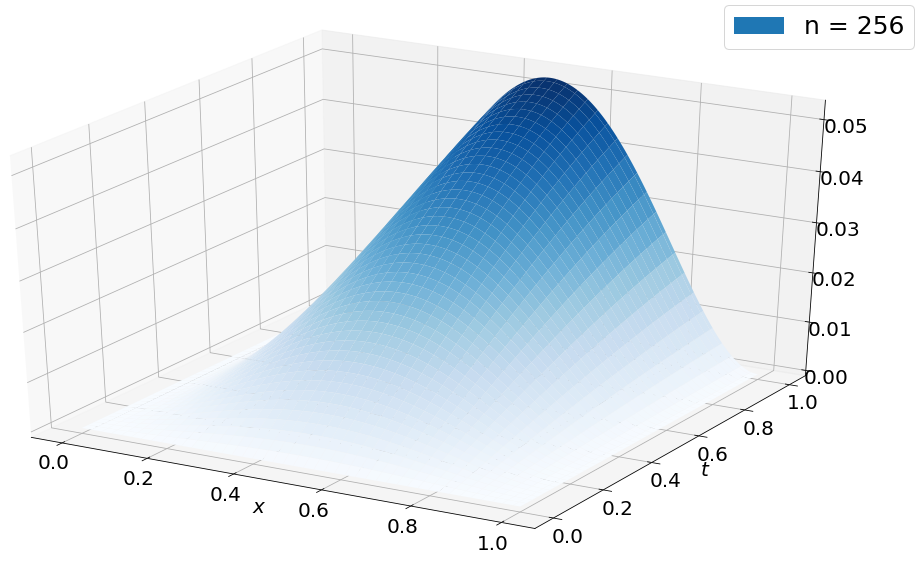}
    }
  \caption{The exact solution (a) with the numerical solution \(u_{\text{num}}\) for \(n=m = 256\) in (b).}
  \label{fig:regular256}
\end{figure}
Thus, let us consider the unique weak solution given by $$u(x,t) = \sin(t) x^2 (x-1)^2,~~(x,t)\in(0,1)\times (0,1),$$ for Problem 1, where \(f_1 = f_2 = 0\), \(c(x) = 1\), \(b(x,t) =1\) and $$g(x,t) = 24\sin(t) + \sin(t)(12x^2 -12x +2)-\sin(t)x^2(x-1)^2.$$

We apply the numerical procedure described above to find the numerical solution \(u_{\text{num}}\) where we set $n=m \in\{32,64,128,256\}$. In Figure \ref{fig:regular256} we plot the exact solution in Figure 1 (a) together with the approximate solution, in Figure 1 (b), obtained by the employed numerical method. 
\begin{figure}[H]
    \centering
    \includegraphics[height=4.5cm]{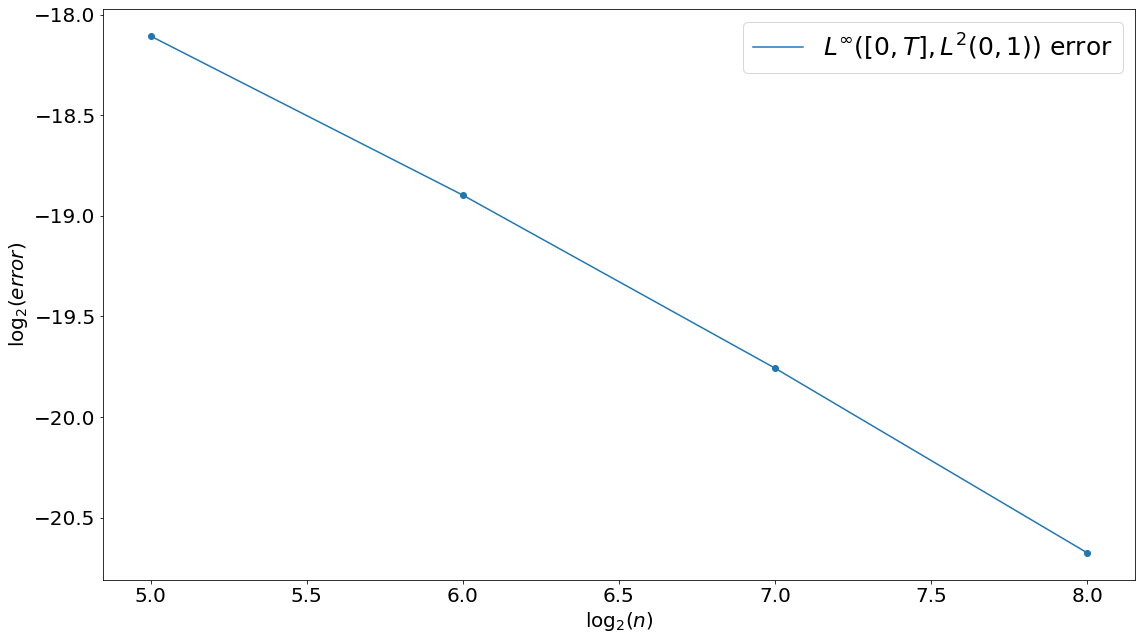}
    \caption{Convergence rate on a logarithmic scale for $n\in\{32,64,128,256\}$.}
    \label{Convergence }
\end{figure}
Moreover, in Figure 2 we plot the convergence rate on a logarithmic scale. We can conclude that the convergence rate for $E_{L^2}$ is approximately of order $3/4$.

\subsection{Irregular cases}

When working with the irregular coefficients, the theoretical guarantees provide us with the existence of the very weak solution. As the coefficients are not smooth enough, we don't have the existence of a classical or weak solution and the only way to obtain the better idea and visualization of the solution is by using the numerical methods. Just as we explained above, the same numerical methods will be employed to plot the numerical solutions and to see that the obtained results are in agreement with the theoretical expectations. We discuss different cases in the following subsections.

\subsubsection{Irregular bending stiffness}

Based on the theoretical results, the very weak solution theorems allow for the bending stiffness coefficient \(c\) to be a positive Radon measure, bounded from below away from zero. We will investigate two irregular cases. First, we consider the case where \(c\) is an \(L^1(0,1)\) function and second, we consider the effect of a Dirac delta term in the bending stiffness. We set the initial conditions equal zero and the other coefficients equal to \(1\). More precisely, \(f_1 = f_2 = 0\), \(b(x,t) = 1\) and \(g(x,t) = 1\).\\

\begin{figure}[H]
  \centering
\subfigure[$\varepsilon=0.01$]{
  \includegraphics[width=.22\textwidth]{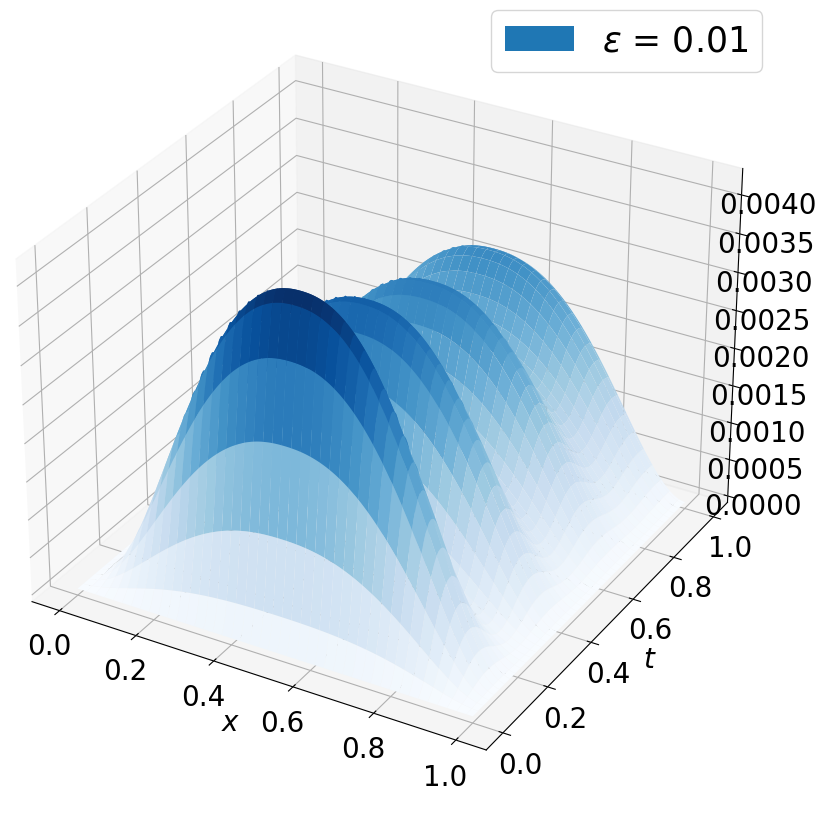}
  }
  \subfigure[$\varepsilon=0.05$]{
  \includegraphics[width=.22\textwidth]{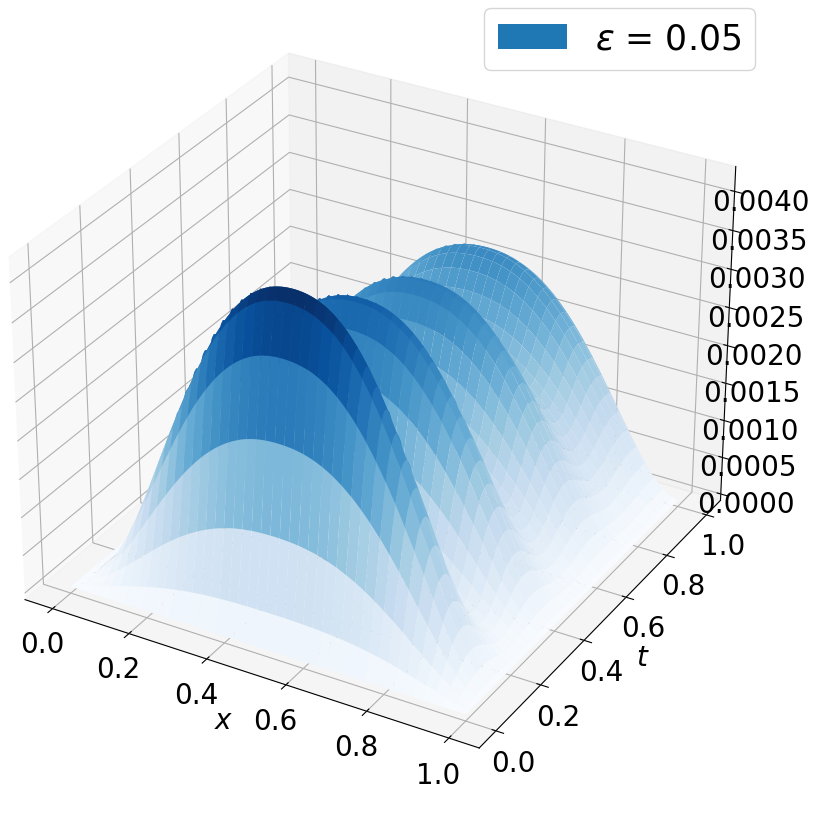}
    }
    \subfigure[$\varepsilon = 0.1$]{
  \includegraphics[width=.22\textwidth]{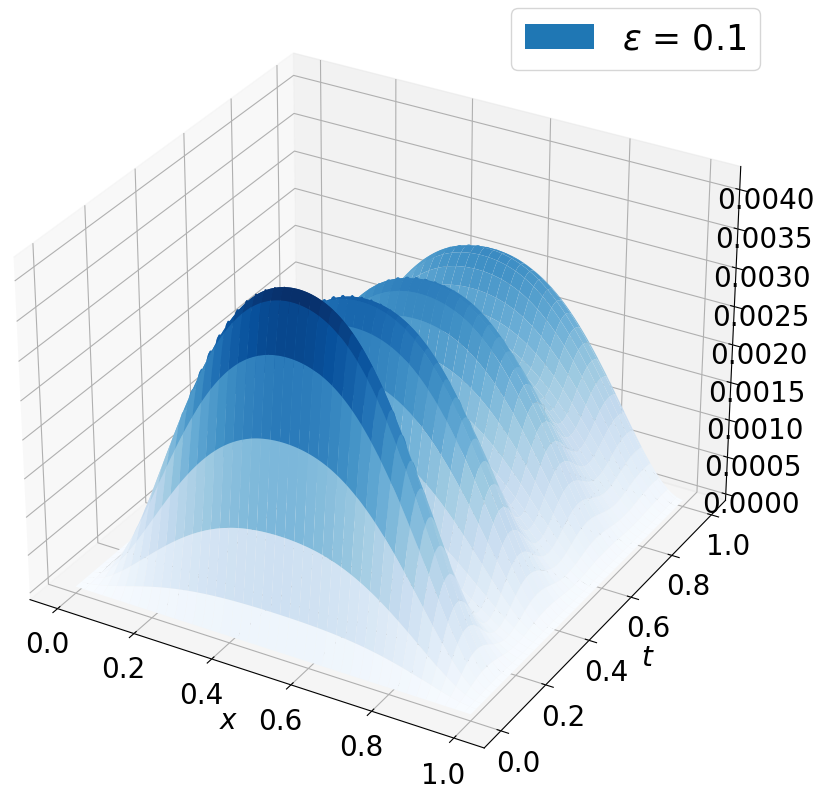}
    }
\subfigure[$\varepsilon = 0.2$]{
  \includegraphics[width=.22\textwidth]{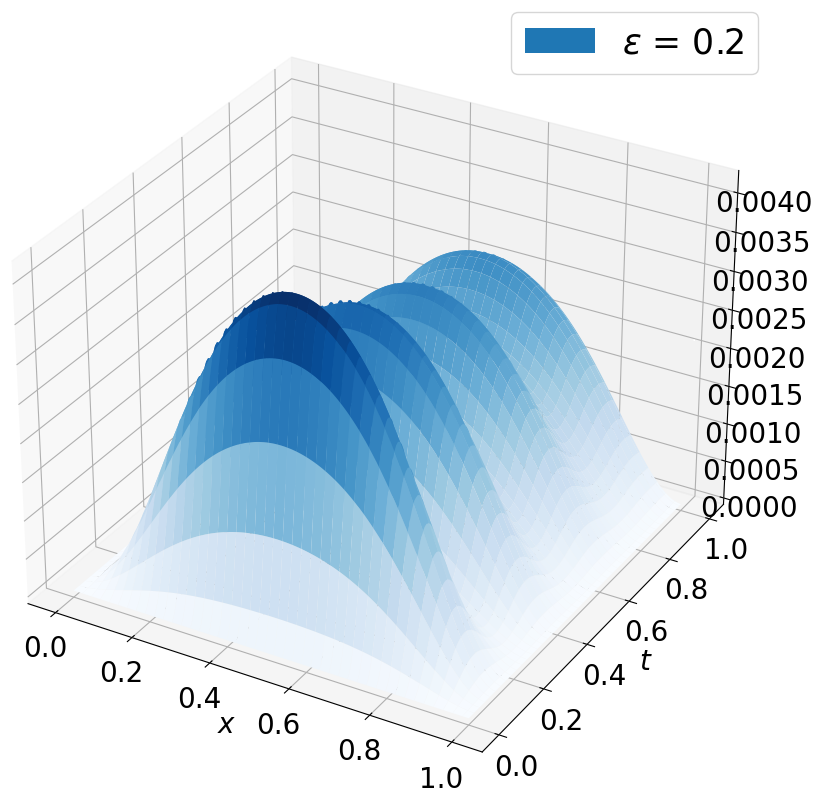}
    }
  \caption{Numerical solution $u_{\varepsilon}(x,t)$ for different values of $\varepsilon$, with $n=256$ and $m=128.$}
  \label{fig:irrlogC}
\end{figure}

1) Let \(c(x) = -\log(|x-0.5|)\), then \(c(x)\) is an \(L^p(0,1)\) function for any \(p\geq 1\), which is not in \(L^\infty(0,1)\). Therefore Theorem \ref{WeakSolution} does not guarantee the existence of a weak solution. However, all the conditions for the existence of a very weak solution are satisfied. In Figure \ref{fig:irrlogC} (a)-(d) we plot the numerical solution $u_{\varepsilon}(x,t)$ for different values of $\varepsilon$, with $c_{\varepsilon}(x)=(\varphi_{\varepsilon}*c)(x)$.
\begin{figure}[H]
    \centering
    \includegraphics[height=4.5cm]{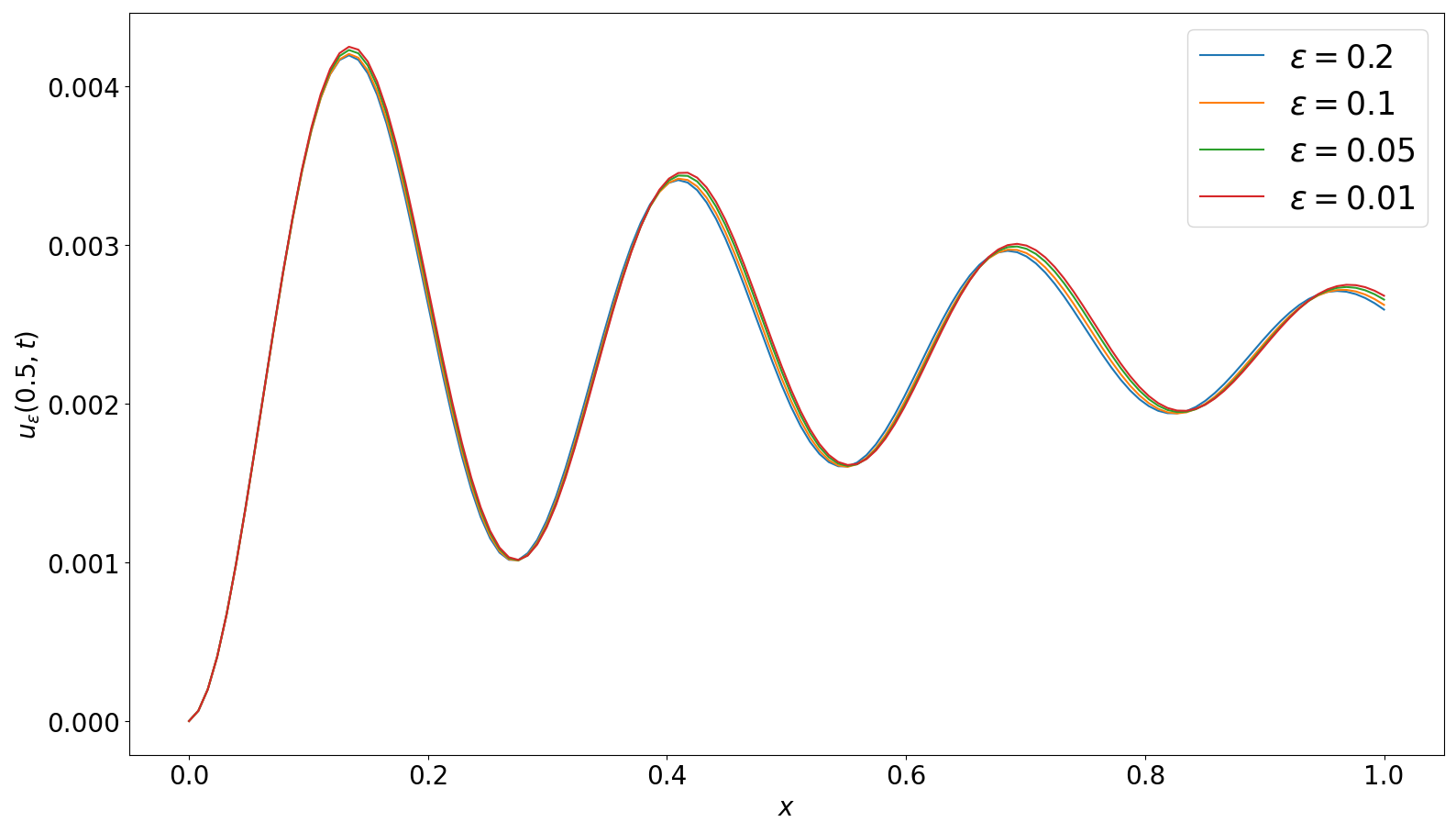}
    \caption{Cross-section for a fixed spatial position $x=0.5$ and $t\in(0,1).$}
    \label{fig:crosslogC}
\end{figure}
Moreover, in Figure \ref{fig:crosslogC} we visualize the cross-section obtained for a fixed spatial position $x=0.5$ and every $t\in(0,1)$. Although a weak solution cannot be guaranteed, we can conclude that as $\varepsilon\to 0$, the solution $u_{\varepsilon}(0.5,t)$ seems to converge.

\begin{figure}[H]
  \centering
\subfigure[$\varepsilon=0.01$]{
  \includegraphics[width=.22\textwidth]{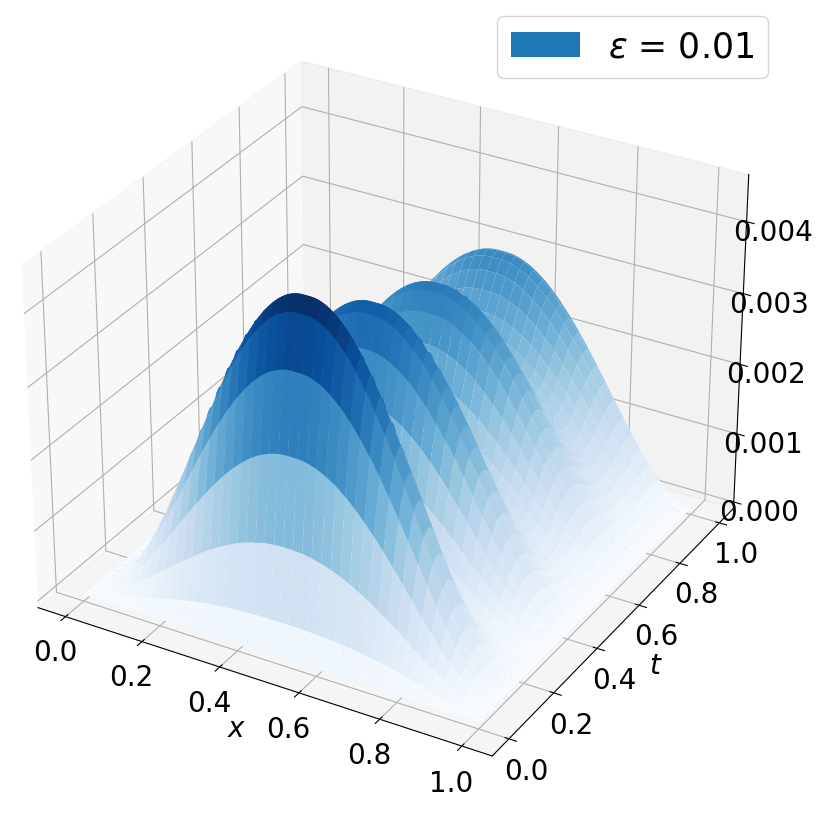}
  }
  \subfigure[$\varepsilon=0.05$]{
  \includegraphics[width=.22\textwidth]{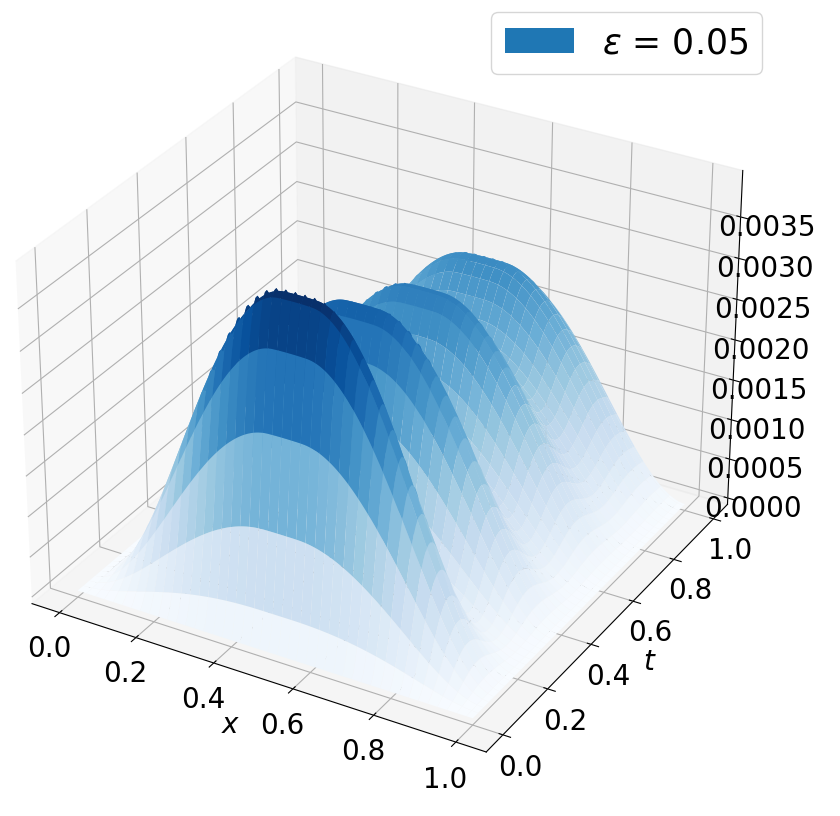}
    }
    \subfigure[$\varepsilon = 0.1$]{
  \includegraphics[width=.22\textwidth]{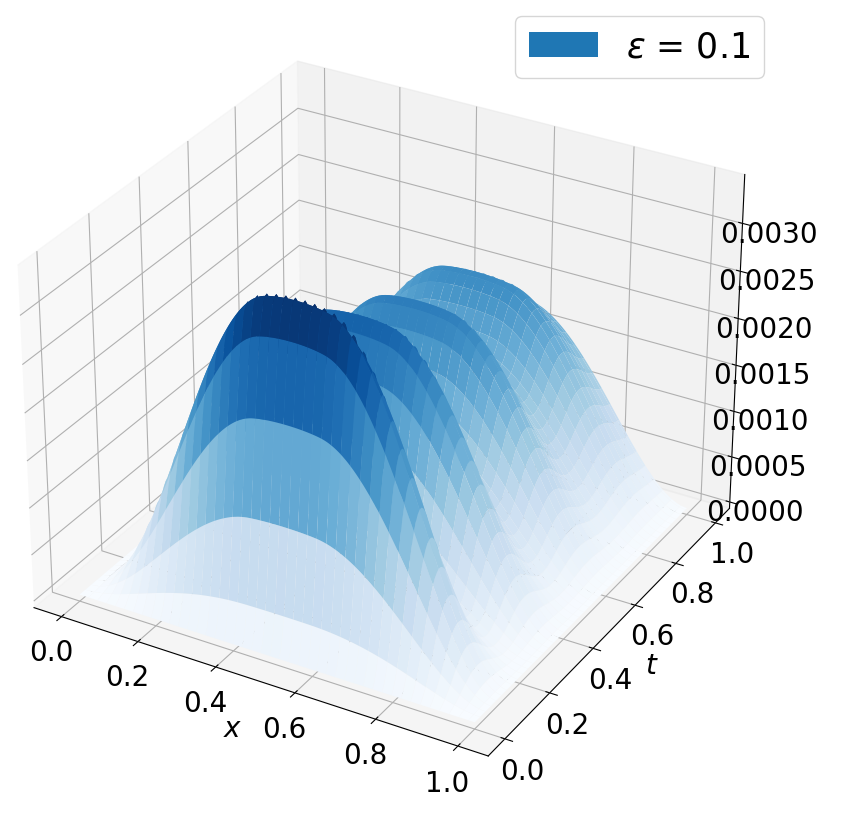}
    }
\subfigure[$\varepsilon = 0.2$]{
  \includegraphics[width=.22\textwidth]{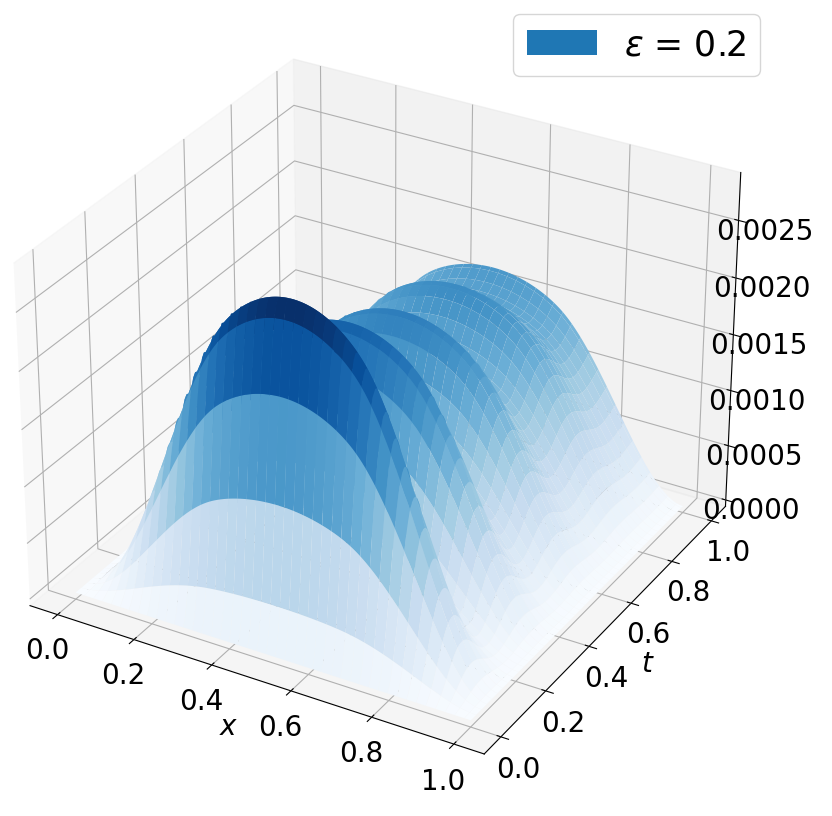}
    }
  \caption{Numerical solution $u_{\varepsilon}(x,t)$ for different values of $\varepsilon$, with $n=256$ and $m=128$.}
  \label{fig:irrC}
\end{figure}

2) In the second case, we consider \(c(x) = 1 + \delta(x-0.5)\) with the regularization given by \(c_\eps(x) = 1+\varphi_\eps(x-0.5)\). Similarly to the previous case, we plot $u_{\varepsilon}(x,t)$ for different values of $\varepsilon$. Those results are visualized in Figure \ref{fig:irrC}.

\begin{figure}[H]
    \centering
   \subfigure[]{
   \includegraphics[width=.45\textwidth]{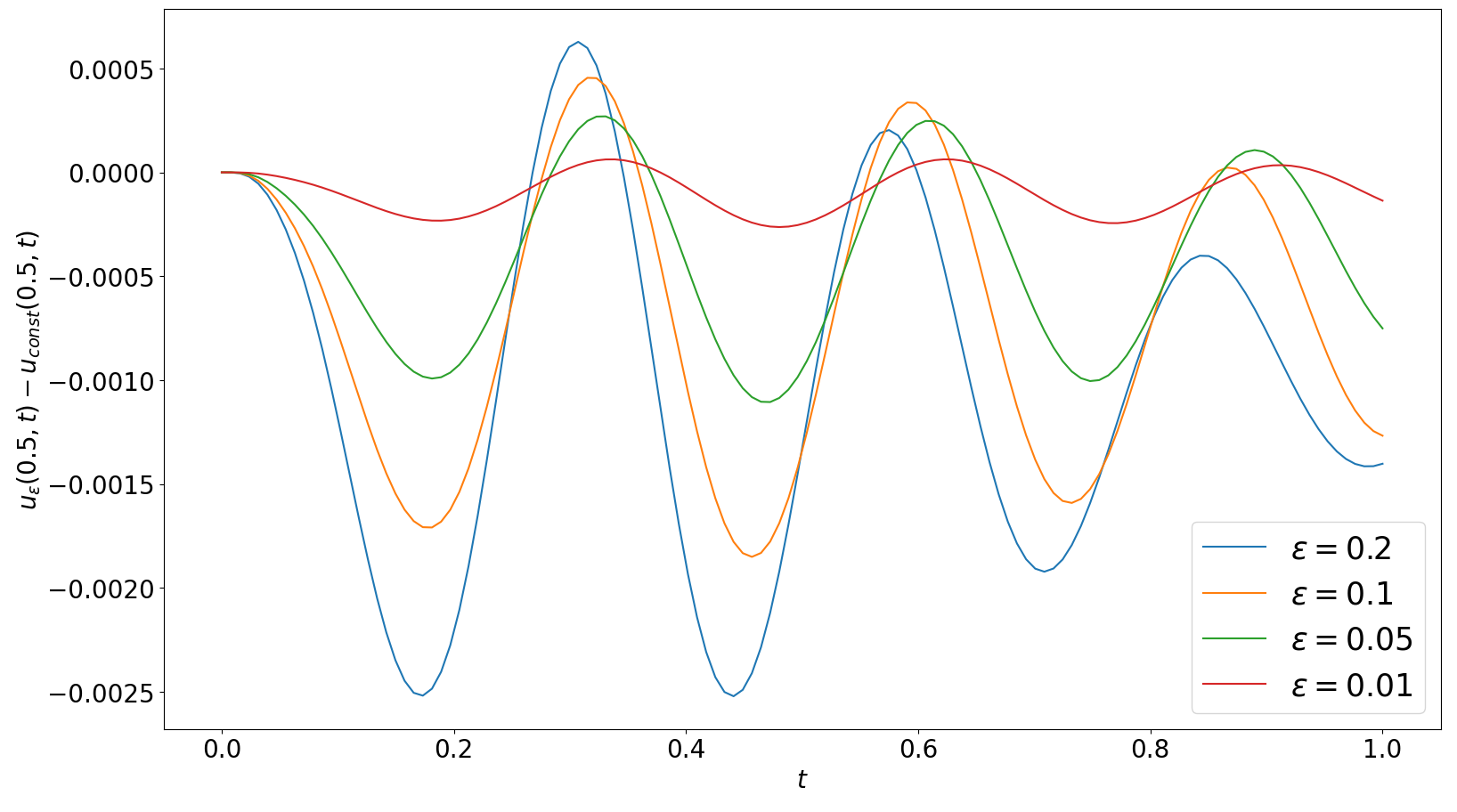}
   }
   \subfigure[]{
   \includegraphics[width=.45\textwidth]{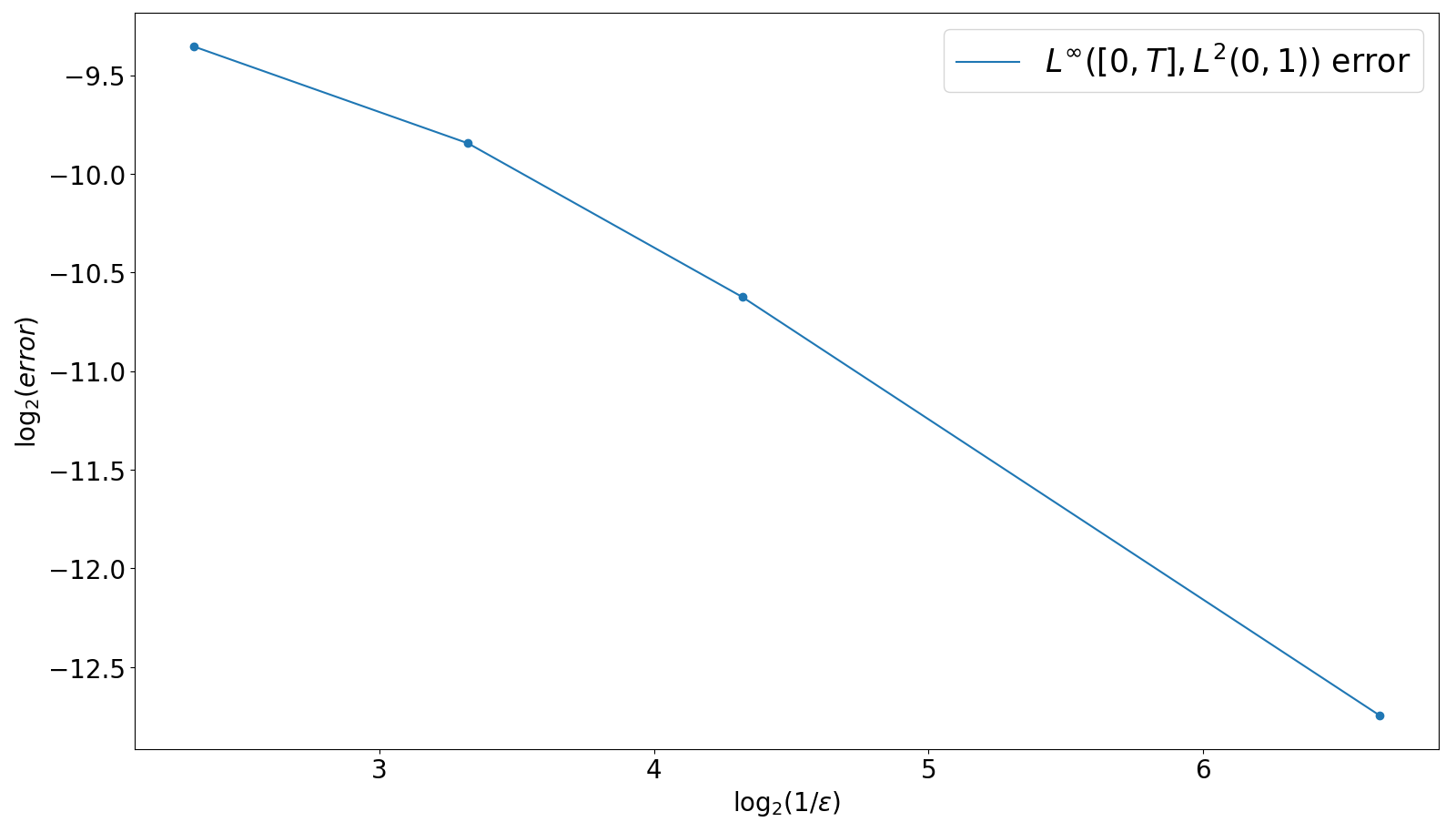}
    }
    \caption{Plot of the cross-section of the difference between $u_{\varepsilon}(0.5,t)$ and $u_{const}(0.5,t)$ in (a), for a fixed spatial position $x=0.5$. In (b) we plot the error $E_{L^2}(u_{\varepsilon},u_{const}$) on a logarithmic scale.}
    \label{fig:errorirrC}
\end{figure}

We propose that the delta term in the bending stiffness has no significant effect on the beam solution. Indeed, consider the case where \(c(x) = 1\) and let us denote the corresponding solution by \(u_{const}\), obtained numerically by the previously discussed procedure. We suspect that the sequential solution converges to $u_{const}$ as $\varepsilon\to 0$, by plotting the cross-section of the difference $u_{\varepsilon}(0.5,t)-u_{const}(0.5,t)$ for $t\in(0,1)$ in Figure \ref{fig:errorirrC} (a), but we confirm this by plotting the error defined in Equation \eqref{errorL2}, between $u_{\varepsilon}$ and $u_{const}$ in Figure \ref{fig:errorirrC} (b). Indeed, based on Fig. 6 (a) we can see that as the $\varepsilon\to 0$, the difference between $u_{\varepsilon}$ and $u_{const}$ decreases. Moreover, in Fig. 6 (b) we can observe that as the $\varepsilon$ decreases, the error tends to zero. 

\subsubsection{Irregular axial force}

\begin{figure}[H]
\centering
\subfigure[$\varepsilon=0.01$]{
  \includegraphics[width=.22\textwidth]{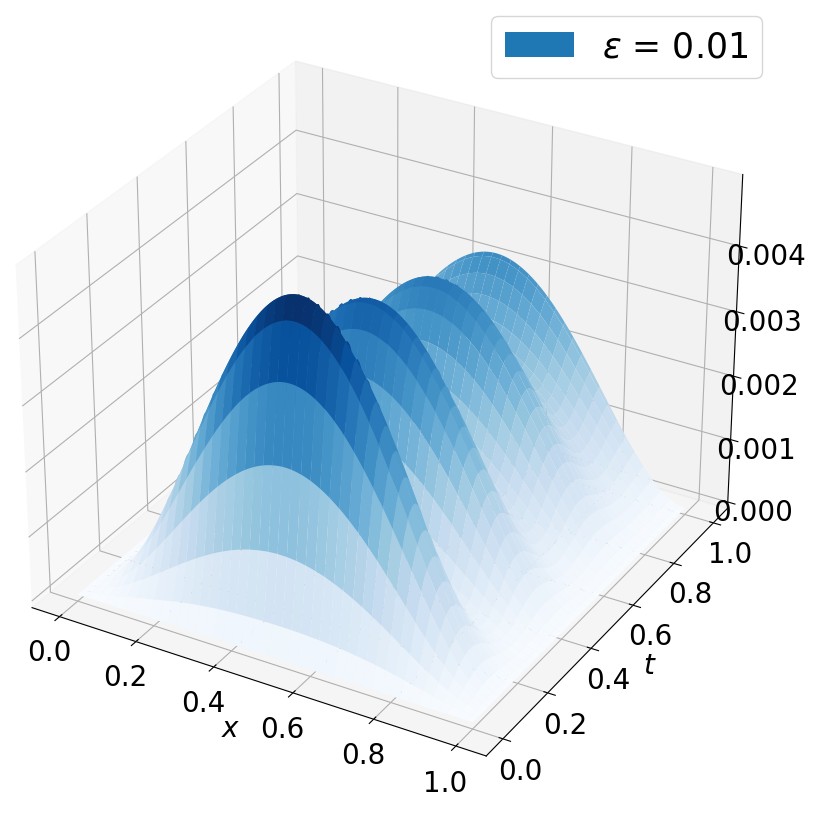}
  }
  \subfigure[$\varepsilon=0.05$]{
  \includegraphics[width=.22\textwidth]{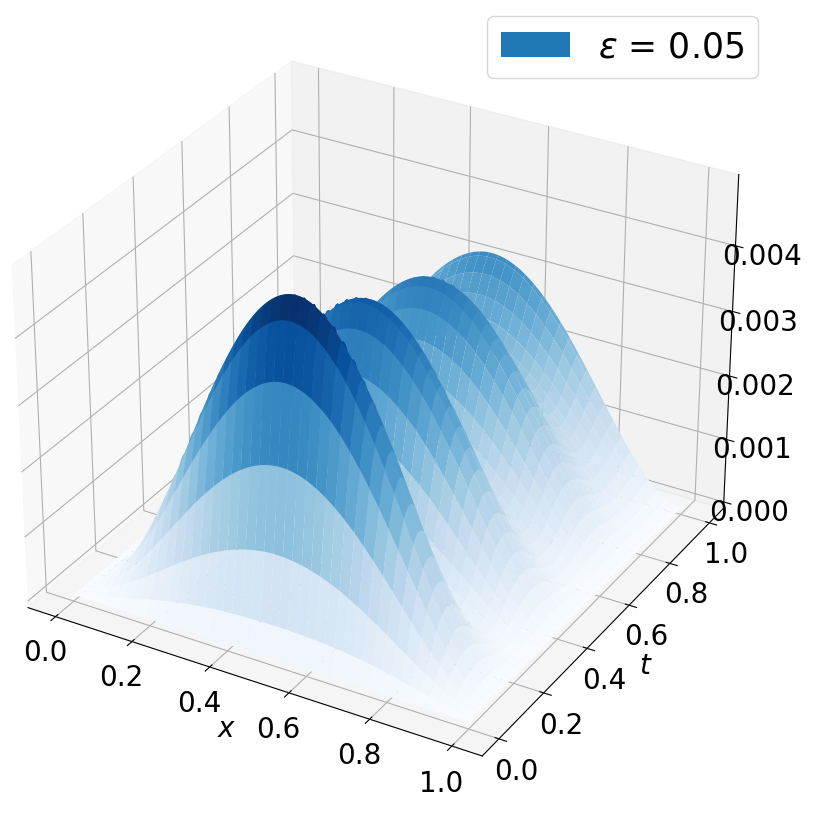}
    }
    \subfigure[$\varepsilon = 0.1$]{
  \includegraphics[width=.22\textwidth]{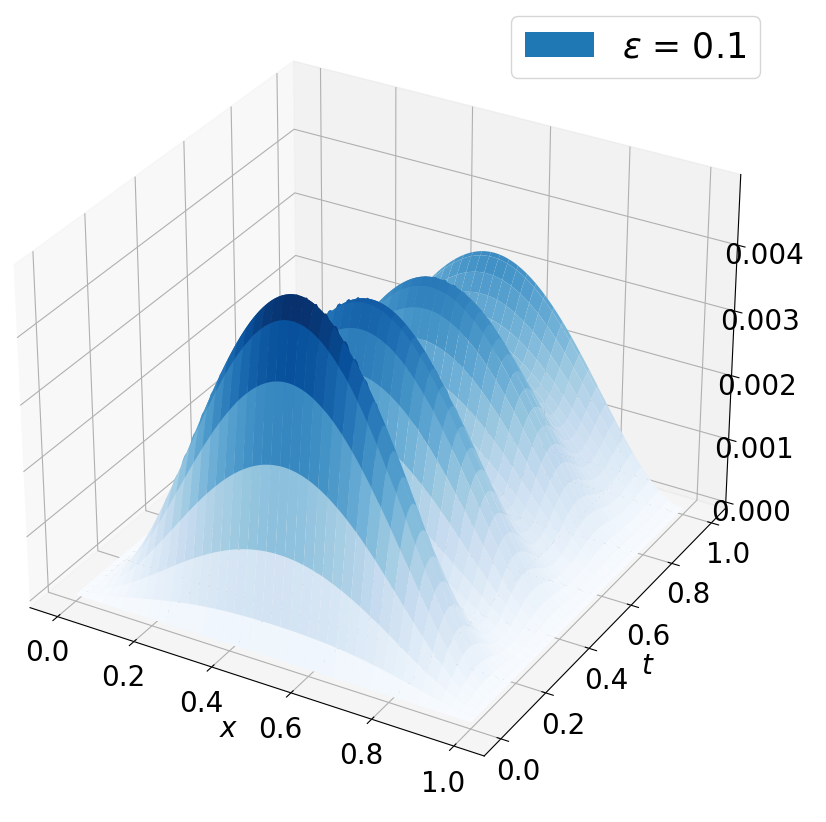}
    }
\subfigure[$\varepsilon = 0.2$]{
  \includegraphics[width=.22\textwidth]{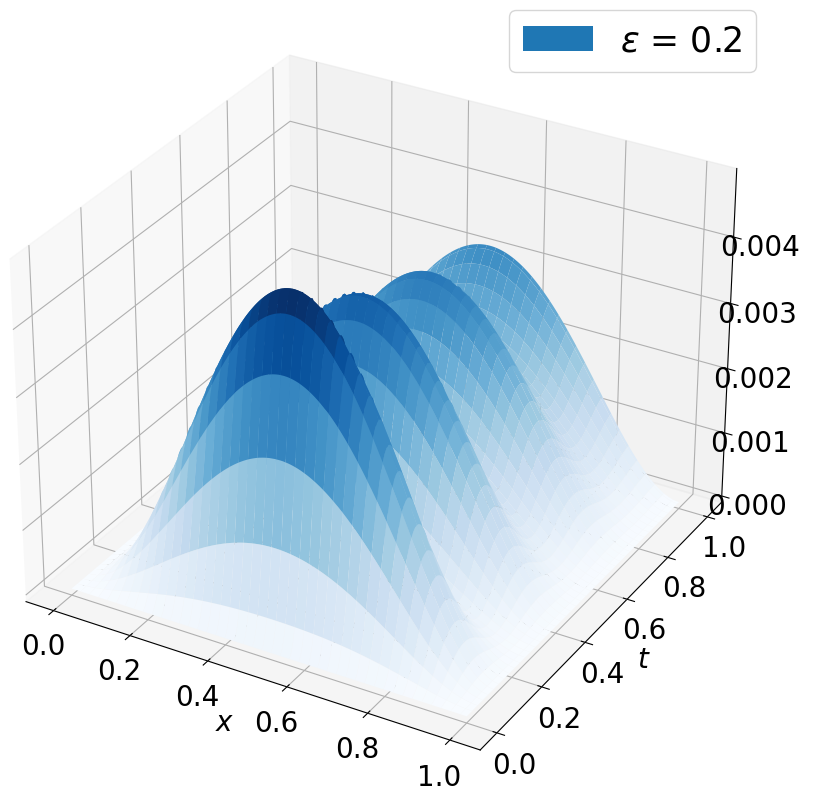}
    }
  \caption{Numerical solution $u_{\varepsilon}(x,t)$ for different values of $\varepsilon$, with $n=256$ and $m=128$.}
  \label{fig:irrB}
\end{figure}

In this subsection, we will investigate an axial force coefficient that consists of a Dirac delta term. This term commonly models an accumulation of axial force at a point due to a crack in the beam material. Thus, let \(b(x,t) = \delta(x-0.5)\) and let us consider the regularization \(b_\eps(x,t) = \varphi_{\varepsilon}(x-0.5)\). 
Again we start by visual investigation in Figure \ref{fig:irrB}, by plotting $u_{\varepsilon}(x,t)$ for different values of $\varepsilon$. 

\begin{figure}[H]
    \centering
   \includegraphics[height=4.5cm]{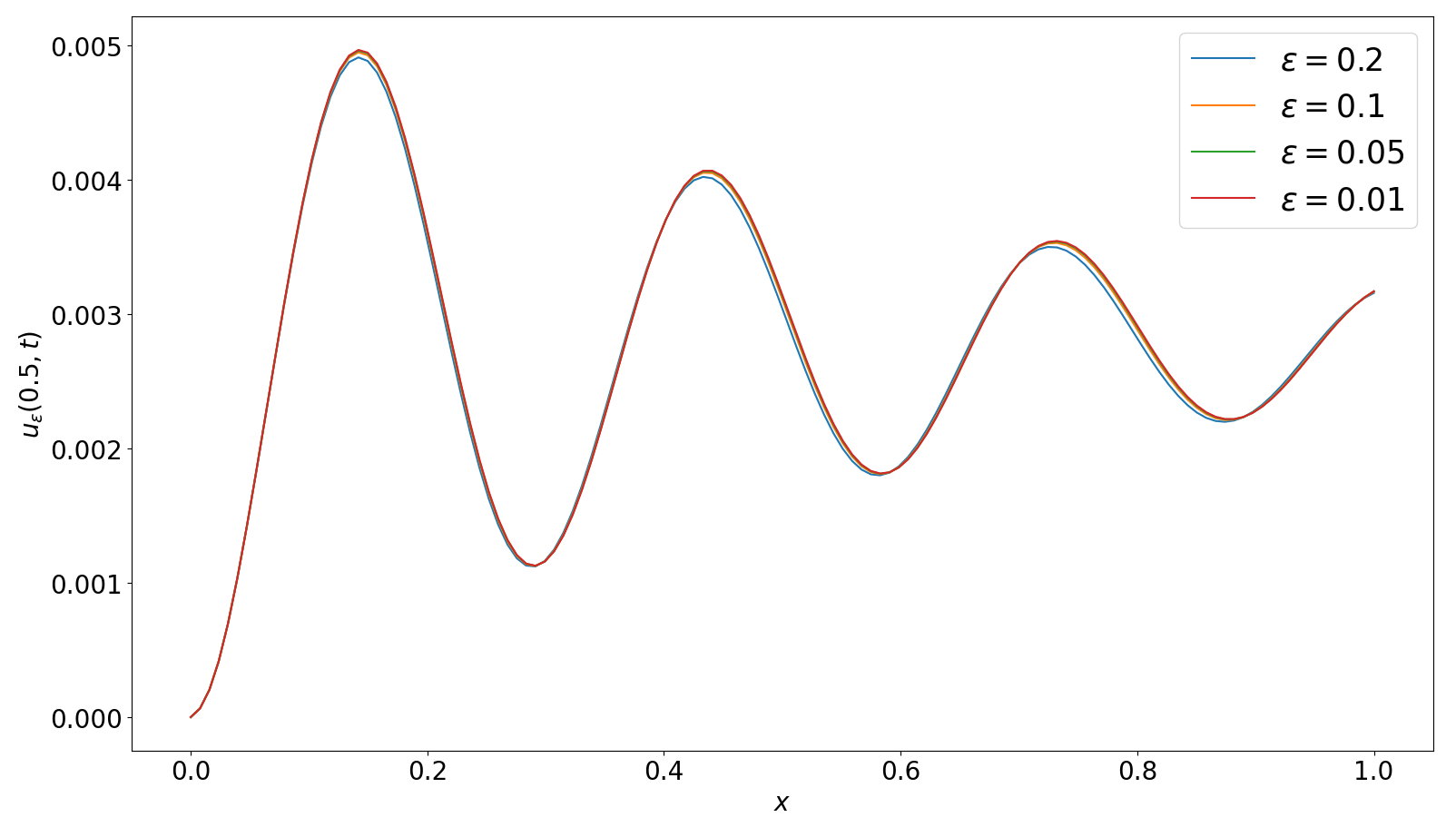}
    \caption{Plot of the cross-section for a fixed spatial position $x=0.5$ and $t\in(0,1).$}
    \label{fig:crossirrB}
\end{figure}

By plotting, in Figure \ref{fig:crossirrB}, we visualize the cross-section obtained for a fixed spatial position $x=0.5$ and every $t\in(0,1)$. Again, the solution $u_{\varepsilon}$ seems to converge. Similar to the previous case of irregular bending stiffness, one can conduct the similar analysis and investigate the effect of the Dirac delta term. 

\subsubsection{Irregular transversal force}

\begin{figure}[H]
  \centering
\subfigure[$\varepsilon=0.01$]{
  \includegraphics[width=.22\textwidth]{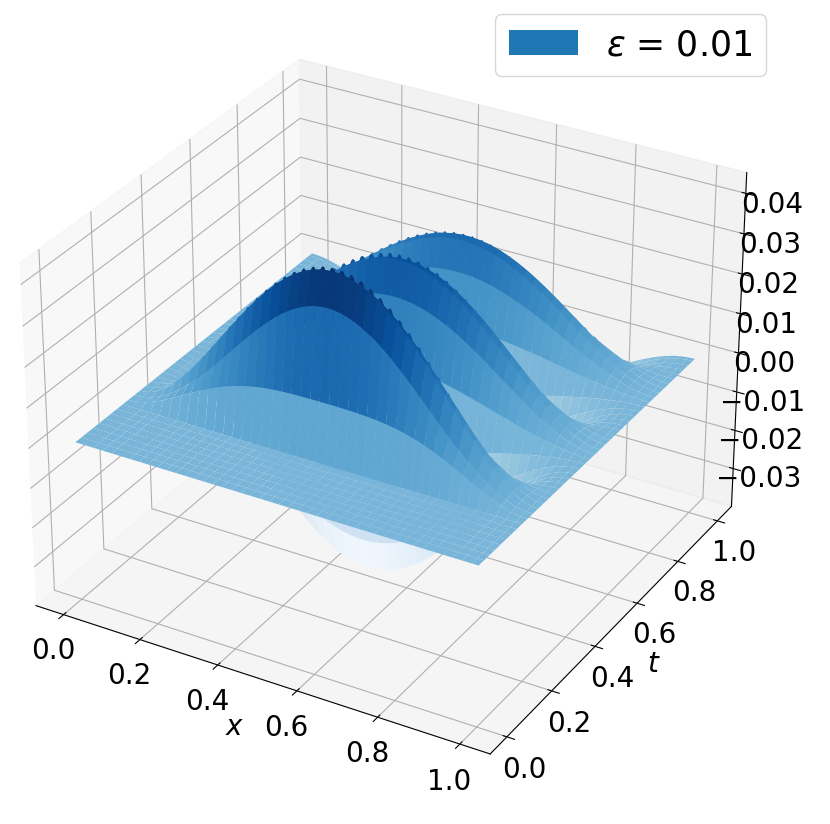}
  }
  \subfigure[$\varepsilon=0.05$]{
  \includegraphics[width=.22\textwidth]{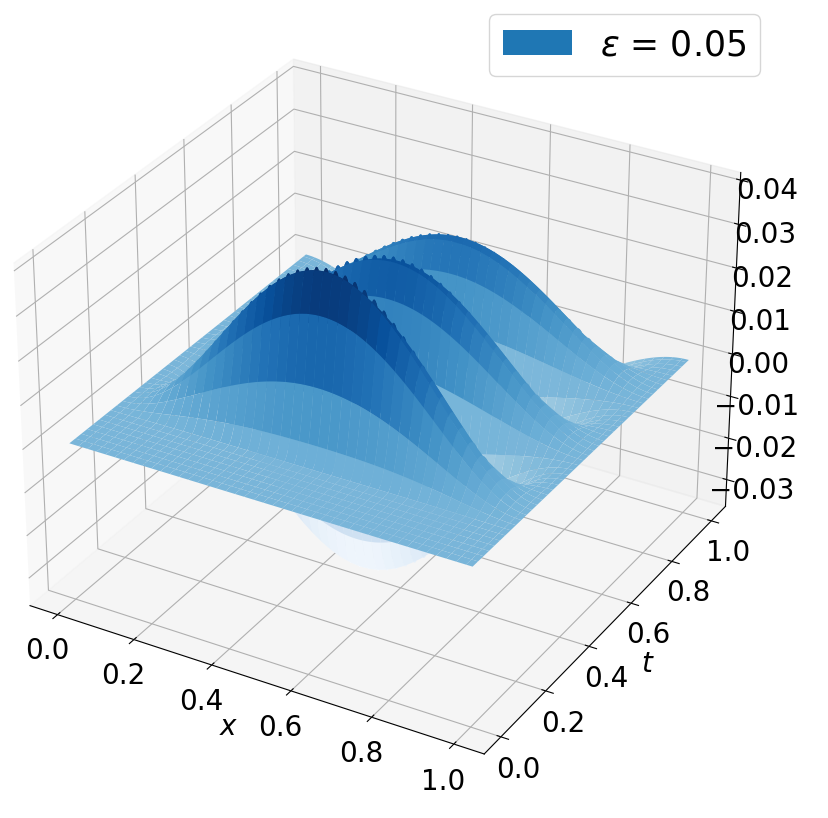}
    }
    \subfigure[$\varepsilon = 0.1$]{
  \includegraphics[width=.22\textwidth]{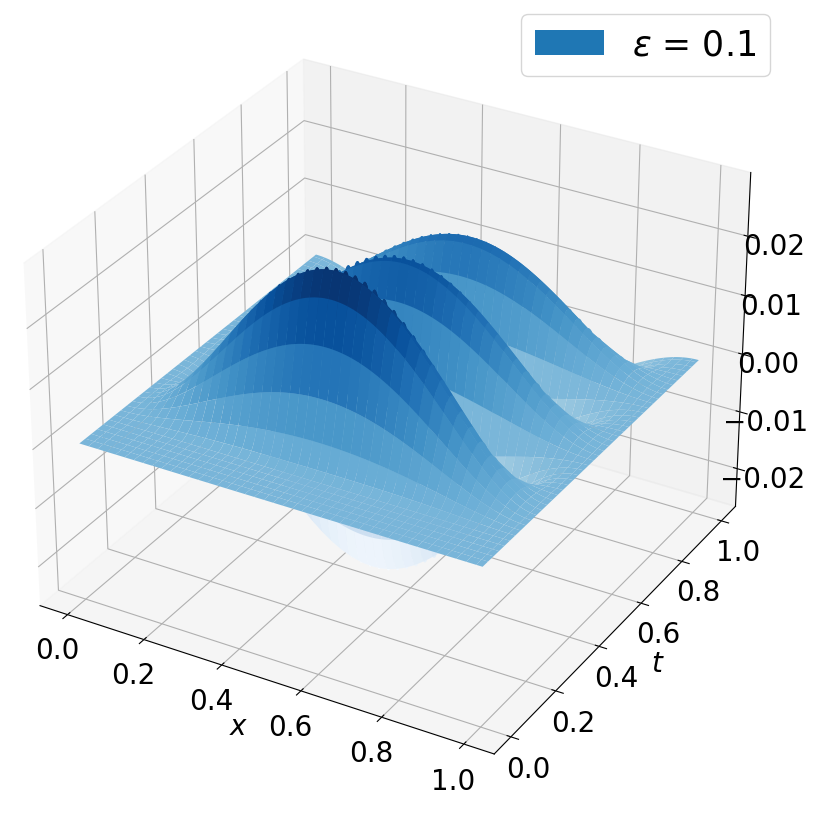}
    }
\subfigure[$\varepsilon = 0.2$]{
  \includegraphics[width=.22\textwidth]{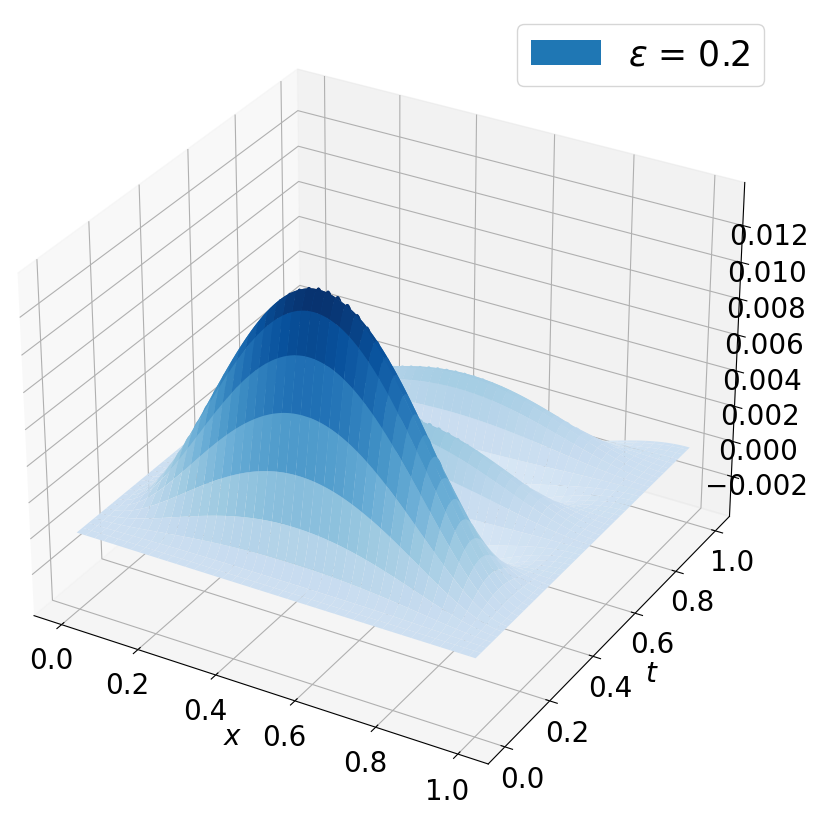}
    }
  \caption{Numerical solution $u_{\varepsilon}(x,t)$ for different values of $\varepsilon$, with $n=128$ and $m=256$.}
  \label{fig:irrg}
\end{figure}

Finally, we discuss the irregular transversal force $g$. First, let us set the initial conditions equal to zero and other coefficients equal to 1. More precisely, let \(f_1 = f_2 = 0\), \(c(x) = 1\) and \(b(x,t) = 1\). Let \(g(x,t) = \delta(t-0.2)\) and consider the regularization \(g_\eps(x,t) = \varphi_\eps(t-0.2)\). Compared to the previous cases, now we set $n=128$ and $m=256$, to obtain more accurate approximation of the Dirac delta term. In Figure \ref{fig:irrg} we visualize the solution $u_{\varepsilon}$ for different $\varepsilon$-values. Compared to the previous cases, we can observe that the solution for $\varepsilon=0.2$ significantly changes compared to the other $\varepsilon$-values.  

\begin{figure}[H]
  \centering
  \includegraphics[height=4.5cm]{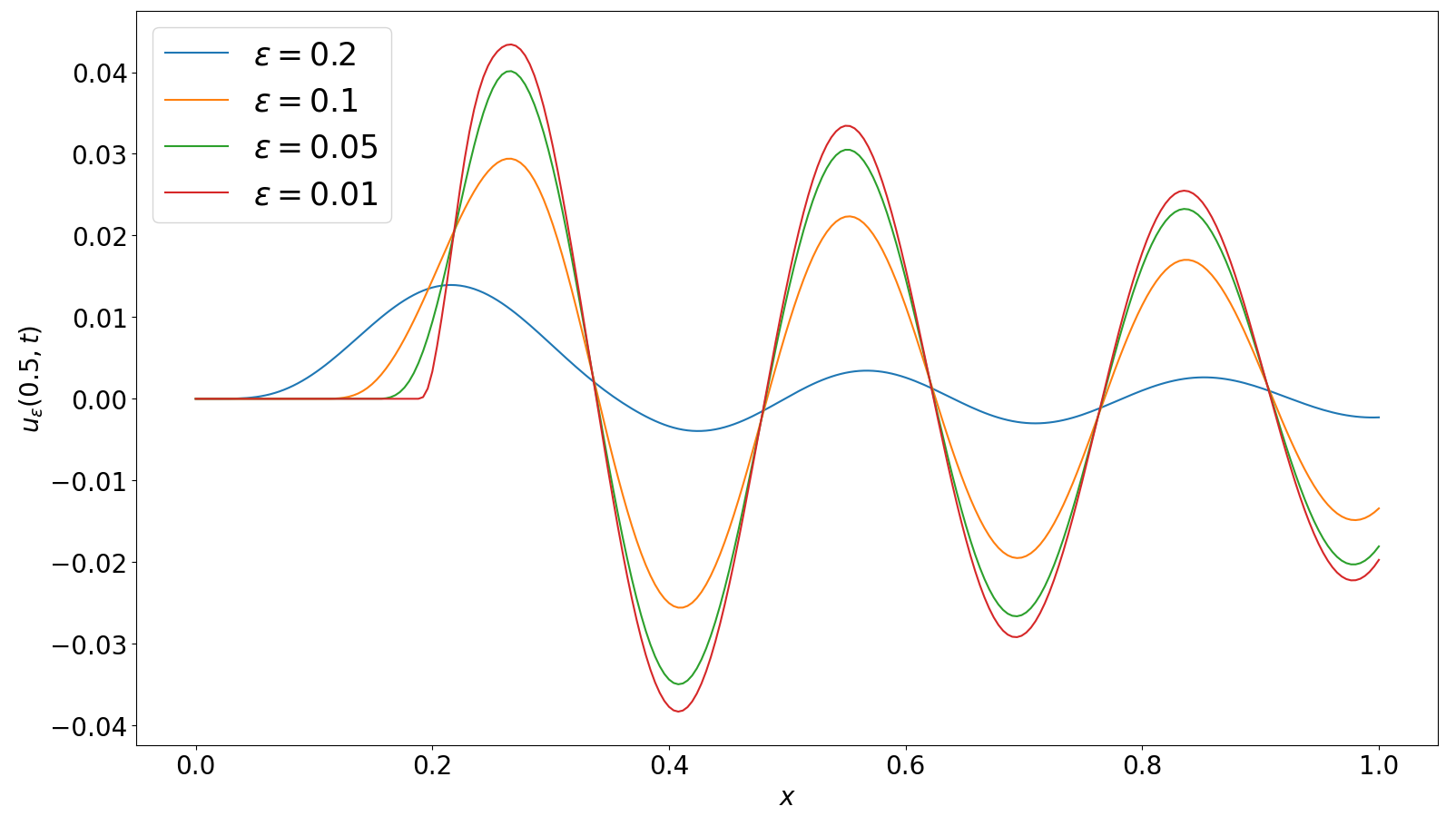}
     \caption{Plot of the cross-section for a fixed spatial position $x=0.5$ and $t\in(0,1).$}  \label{fig:deltaTGCrosssection}
\end{figure}

Taking the cross-section for a fixed spatial position $x=0.5$ and $t\in (0,1)$ we can clearly see that the delta term forces a peak at $t=0.2$, which becomes more pronounced as $\varepsilon$ goes to zero. Those findings are visualized in Figure \ref{fig:deltaTGCrosssection}.

\section{Conclusion}
\label{concl}

In this article we show the existence, uniqueness and consistency results of the very weak solution for the Euler-Bernoulli equation with distributional coefficients and initial conditions. We consider nets of non-smooth data to obtain a sufficiently regular solution net. Our treatment of this problem, broadens the scope of the very weak solutions concepts. For the Euler-Bernoulli problem, we extend the existing weak solution framework making it suitable for the very weak solution approach. We defined and showed existence of the very weak solution, together with the uniqueness results defined in an appropriate sense. Moreover, we give the characterization of the distributions that allow the existence of the solution of the discussed problem. Lastly, we showed the consistency of the very weak solution with the weak solution, when the irregular data were sufficiently nice. Numerical analysis conducted showed that the very weak solution coincides with the weak solution, when the latter exists. When the coefficients are irregular functions, we showed that numerical analysis can provide us with more insights into the behaviour of the  very weak solution.

\end{document}